\documentclass[a4paper,12pt]{amsart}
\usepackage{graphicx}
\usepackage{pdfsync}
\usepackage[leqno]{amsmath}
\usepackage{amsfonts,amssymb}
\usepackage[all]{xy}
\usepackage{pb-diagram,pb-xy}
\numberwithin{equation}{section}

\theoremstyle{plain}
\newtheorem{prop}{Proposition}[section]
\newtheorem{lem}[prop]{Lemma}
\newtheorem{thm}[prop]{Theorem}
\newtheorem{cor}[prop]{Corollary}

\newtheorem*{thmapp}{Theorem \ref{catso10=21}}

\theoremstyle{definition}
\newtheorem{defi}[prop]{Definition}

\theoremstyle{remark}
\newtheorem{rem}{Remark}

\newcommand{\cat}[1]{\operatorname{cat}(#1)}
\newcommand{\cupln}[2]{\operatorname{cup}(#1;#2)}
\newcommand{\cl}[1]{\operatorname{Cat}(#1)}
\newcommand{\ad}[1]{\operatorname{ad}(#1)}
\newcommand{\id}[1]{\operatorname{id}_{#1}}
\newcommand{\proj}{\operatorname{proj}}
\newcommand{\incl}{\operatorname{incl}}
\def\field{\mathbb{F}}
\def\hookxrightarrow#1{\lhook\joinrel\xrightarrow{#1}}

\begin{document}
\title{On Lusternik-Schnirelmann category of $\mathrm{SO}(10)$}
\author[Iwase]{Norio Iwase${}^{\dagger}$}
\author[Kikuchi]{Kai Kikuchi}
\author[Miyauchi]{Toshiyuki Miyauchi}
\email[Iwase]{iwase@math.kyushu-u.ac.jp}
\thanks{${}^{\dagger}$supported by the Grant-in-Aid for Scientific Research \#22340014 from Japan Society for the Promotion of Science.}
\email[Miyauchi]{miyauchi@math.sci.fukuoka-u.ac.jp}
\address[Iwase]{Faculty of Mathematics,
Kyushu University,
Motooka 744, Fukuoka 819-0395, Japan}
\address[Miyauchi]{Department of Applied Mathematics,
Faculty of Science,
Fukuoka University,
Fukuoka, 814-0180, Japan}
\date{\today}
\begin{abstract}
Let $G$ be a compact connected Lie group and $p : E\to \Sigma A$ be a principal G-bundle with a characteristic map $\alpha : A\to G$, where $A=\Sigma A_{0}$ for some $A_{0}$.
Let $\{K_{i}{\to} F_{i-1}{\hookrightarrow} F_{i} \,|\, 1{\le} i {\le} n,\, F_{0}{=} \{\ast\} \; F_{1}{=} \Sigma{K_{1}} \; \text{and}\; F_{n}{\simeq} G \}$ be a cone-decomposition of $G$ of length $m$ and $F'_{1}=\Sigma{K'_{1}} \subset F_{1}$ with $K'_{1} \subset K_{1}$ which satisfy $F_{i}F'_{1} \subset F_{i+1}$ up to homotopy for any $i$.
Our main result is as follows: we have $\cat{X} \le m{+}1$, if firstly the characteristic map $\alpha$ is compressible into $F'_{1}$, secondly the Berstein-Hilton Hopf invariant $H_{1}(\alpha)$ vanishes in $[A, \Omega F'_1{\ast}\Omega F'_1]$ and thirdly $K_{m}$ is a sphere.
We apply this to the principal bundle $\mathrm{SO}(9)\hookrightarrow\mathrm{SO}(10)\to S^{9}$ to determine L-S category of $\mathrm{SO}(10)$.
\end{abstract}
\maketitle
\section{Introduction}
\label{intro}
In this paper, we work in the homotopy category of pointed $CW$-complex.
So we identify a map with its homotopy class, unless it causes a confusion.
The Lusternik-Schnirelmann category of a space $X$ is the least integer $n$ such that there exists an open covering 
$U_{0},\dots, U_{n}$ of $X$ with each $U_{i}$ contractible in $X$ and is denoted by $\cat{X}$.
If no such integer exists, we write $\cat{X}=\infty$.
\begin{thm}[Ganea \cite{Ganea67}]
\label{thm:ganea}
Let $X$ be a connected space.
Then there is a sequence of fibrations $F_{n}X \hookrightarrow G_{n}X \to X$, natural with respect to $X$ so that $\cat{X} \leq n$ if and only if the fibration $G_{n}X \to X$ has a cross-section.
\end{thm}
We can replace the inclusion $F_{n}X \hookrightarrow G_{n}X$ by the fibration $p_{n}^{\Omega{X}} : E^{n+1}\Omega{X} \to P^{n}\Omega{X}$ associated with the $A_{\infty}$ structure of $\Omega{X}$ in the sense of Stasheff \cite{Stasheff63}, where $E^{n+1}\Omega X$ has the homotopy type of $\Omega X^{\ast(n+1)}$ the $n{+}1$-fold join of $\Omega X$ and $P^{n}\Omega{X}$ is the $\Omega X$-projective $n$-space.
We know that there is a natural equivalence $P^{\infty}\Omega{X} \simeq X$ together with $P^{0}\Omega{X}=\ast$ and $P^{1}\Omega{X}=\Sigma\Omega{X}$.
$P^{n}\Omega{X}$ is also equipped with a map $e_n^X$ given by the composition $P^n\Omega X \hookrightarrow P^{\infty}\Omega X \simeq X$, which is natural at $X$: for any map $f : X \to Y$, we have $f \circ e^{X}_{n} = e^{Y}_{n} \circ P^{n}\Omega f$.
We know that $e^X_1$ is the evaluation map, i.e, $e^{X}_{1}=ev : \Sigma\Omega{X} = P^{1}\Omega{X} \to X$ given by $ev(t{\wedge}\ell)=\ell(t)$.
Theorem \ref{thm:ganea} can be reformulated in the following form.
\begin{thm}
\label{thm-ganea67}
$\cat{X} \leq n$ if and only if $e^{X}_{n} : P^{n}\Omega{X} \to X$ has a right homotopy inverse.
\end{thm}
\par
Let $R$ be a commutative ring and $X$ a connected space.
The cup-length of $X$ with coefficients in $R$ is the least non-negative integer $k$ (or $\infty$) such that all ($k{+}1$)-fold cup products vanish in the reduced cohomology $\Tilde{H}^{\ast}(X;R)$.
We denote this integer $k$ by $\cupln{X}{R}$ following Iwase \cite{Iwase04}.
\par
In 1967, Ganea introduced in \cite{Ganea67} a homotopy invariant $\cl{X}$ for a space $X$, modifying Fox's strong category.
In the same paper, he gave the following characterization using the notion of a cone-decomposition.
\begin{defi}[Ganea \cite{Ganea67}]
The strong category $\cl{X}$ of a connected space $X$ is $0$ if $X$ is contractible and, otherwise, is equal to the least integer $n$ such that there is a sequence of cofibrations (a cone-decomposition of length $m$)
\[
\{K_{i}\to F_{i-1}\hookrightarrow F_{i} \;|\; 1 \le i \le m\}
\]
with $F_{0}= \{\ast\}$ and $F_{m}\simeq X$.
$\cl{X}$ is often called the cone-length of $X$.
\end{defi}
The following inequalities among these invariants are well-known
\[ \cupln{X}{R}\leq \cat{X} \leq \cl{X}. \]
\par
Let $f : \Sigma X \to \Sigma Y$ be a map. We denote 
$H_1(f) \in [\Sigma X, \Omega \Sigma Y{\ast} \Omega \Sigma Y]$ 
by the Berstein-Hilton Hopf invariant (see Berstein and Hilton \cite{BH60}).
\par
The purpose of this paper is to prove the following theorems: let $G\hookrightarrow E\rightarrow \Sigma A$, where $A=\Sigma A_0$ for some $A_{0}$, be a principal bundle with a characteristic map $\alpha : A\rightarrow G$. 
We assume that $G$ is a connected compact Lie group with a cone-decomposition of length $m$, that is, there are sequence of cofibrations
\[
\{K_i\rightarrow F_{i-1}\hookrightarrow F_i \,|\,1\leq i\leq m \}
\]
with $F_0=\{\ast\}$, $F_{1} \simeq \Sigma{K_{1}}$ and $F_m \simeq G$, and so we denote by $\mu : F_{m} \times F_{m} \to F_{m} \simeq G$ the multiplication of $G$.
\begin{thm}\label{main2thm}
Let $F^{\prime}_1=\Sigma K^{\prime}_1$ where $K^{\prime}_1$ is a subspace of $K_1$ so that $F'_{1} \subset F_{1}$ and let $\mu\vert_{F_{i}{\times}F'_{1}} : F_{i} \times F'_{1} \to F_{m}$ be a map compressible into a map $\mu'_{i,1} : F_{i} \times F'_{1} \to F_{i+1}$ for any $i$ with $1\leq i < m$ such that $\mu'_{i,1}\vert_{F_{i-1}{\times}F'_{1}} \simeq \mu'_{i-1,1}$ in $F_{i+1}$.
Then the following three conditions imply $\cat{E}\leq m{+}1$.
\begin{enumerate}
\item
$\alpha$ is compressible into $F^{\prime}_1$
\item
$H_1(\alpha)=0$ in $[A, \Omega F'_1{\ast} \Omega F'_1]$,
\item
$K_m$ is a sphere.
\end{enumerate}
\end{thm}
\begin{cor}\label{main1thm}
Let $\mu\vert_{F_{i}{\times}F_{1}} : F_{i} \times F_{1} \to F_{m}$ be a map compressible into a map $\mu_{i,1} : F_{i} \times F_{1} \to F_{i+1}$ for any $i$ with $1\leq i < m$ such that $\mu_{i,1}\vert_{F_{i-1}{\times}F_{1}} \simeq \mu_{i-1,1}$ in $F_{i+1}$.
Then the following three conditions imply $\cat{E}\leq m{+}1$.
\begin{enumerate}
\item
$\alpha$ is compressible into $F_1$,
\item
$H_1(\alpha )=0$ in $[A, \Omega F_1{\ast} \Omega F_1] $,
\item
$K_m$ is a sphere.
\end{enumerate}
\end{cor}
\begin{thmapp}
$\cat{\mathrm{SO}(10)}=21=\cupln{\mathrm{SO}(10)}{\field_{2}}$.
\end{thmapp}
This would suggest that $\cat{\mathrm{SO}(n)}=\cupln{\mathrm{SO}(n)}{\field_{2}}$ for all $n$.

In Sections \ref{structure} and \ref{cone-decom}, 
we construct a structure map and a cone-decomposition of some spaces which play the vital role in the proof of the main theorems.
In Section \ref{str-cone}, we show the crucial relation between a structure map and a cone-decomposition which are constructed in Section \ref{structure} and \ref{cone-decom}.
In Section \ref{proof-mthm}, we prove Theorem \ref{main2thm}.
Finally in Section \ref{appli-mthm}, we determine $\cat{\mathrm{SO}(10)}$.
%
\section{Structure map associated with a filtration}
\label{structure}
\begin{defi}
A space $X$ equipped with a sequence of subspaces $\{X_{n};n{\geq}0\}$,
\[
X \supset\cdots\supset X_n \supset X_{n-1} \supset\cdots\supset \{\ast\}
\]
is called a space $X$ filtered by $\{X_{n};n{\geq}0\}$ or simply a filtered (based) space, and is denoted by $=(X,\{X_{n};n{\geq}0\})$.
We also denote by $i^{X}_{m,n} : X_m\to X_n$, $m < n$ the inclusion map of filtrations.
\end{defi}
\begin{defi}
Suppose that the space $X$ and $Y$ are filtered by $\{ X_n \}$ and $\{ Y_n \}$, respectively.
A filtered map $f : X\to Y$ is a filtration-preserving map, that is,
$f(X_n)\subset Y_n$ for all $n$.
\end{defi}
We denote $p_{m}^{\Omega X}$ by the map $E^m\Omega X \to P^{m-1}\Omega X$
in Theorem \ref{thm-ganea67} and $\iota^{\Omega X}_{m,n} : P^{m}\Omega X\to P^{n}\Omega X$
 by the inclusion map for $m < n$.
\begin{prop}\label{cd-to-gs}
Let $X$ and $Y$ be filtered by $\{ X_n \}$ and $\{ Y_n \}$, respectively, where $\{ X_n \}$ is a cone-decomposition of $X$ by the sequence $\{L_{i} \xrightarrow{h_{i}} X_{i-1} \hookrightarrow X_{i}\,|\,1 {\leq} i {\leq} n\}$ of cofibrations.
If a map $f : X \to Y$ is a filtered map w.r.t. the above filters, then there exist families of maps $\{\hat{f}_{i} : X_{i}\to P^{i}\Omega Y_{i} \,|\, 0 {\leq} i {\leq} n \}$ and $\{\hat{f}^{0}_{i} : L_{i}\to E^{i}\Omega Y_{i} \,|\,1 {\leq} i {\leq} n \}$ such that $\{\hat{f}_{i}\}$ and $\{\hat{f}^{0}_{i}\}$ satisfy the following conditions.
\begin{enumerate}
  \item Let $i_{i-1,i}^X : X_{i-1} \hookrightarrow X_{i}$ be the canonical inclusion.
 Then the following diagram is commutative.                    
$$\xymatrix{
L_{i} \ar[r]^{h_{i}} \ar[dd]_{\hat{f}^{0}_{i}}
 & X_{i-1} \ar@{^{(}->}[r]^{i_{i-1,i}^X} \ar[d]^{\hat{f}_{i-1}}
 & X_{i}  \ar[dd]^{\hat{f}_{i}}
\\
 & P^{i-1}\Omega Y_{i-1} \ar@{^{(}->}[d]^{P^{i-1}\Omega i_{i-1,i}^Y}
 &
\\
E^{i}\Omega Y_{i} \ar[r]_{P_{i}^{\Omega Y_i}}
 & P^{i-1}\Omega Y_{i} \ar@{^{(}->}[r]_{\iota_{i-1,i}^{\Omega Y_i}}
 & P^{i}\Omega Y_{i}.
}$$
 \item $e_{i}^{Y_{i}}{\circ}\hat{f}_{i}=  f|_{X_{i}}$
\end{enumerate}
\end{prop}
\begin{proof}
First of all, we define $\hat{f}_0$ as the trivial map, i.e, $\hat{f}_0=*$.
\par
Next, we proceed by induction on $i\geq1$.
When $i=1$, we define $\hat{f}^{0}_1=\ad{f|_{X_1}}$ and $\hat{f}_1=\Sigma \ad{f|_{X_1}}$, so that the following diagram is commutative.
$$\xymatrix{
L_{1} \ar[r] \ar[d]_{\hat{f}^{0}_{1}}
 & \ast \ar[r] \ar[d]^{\hat{f}_0}
 & \Sigma L_1 \ar[d]^{\hat{f}_{1}}\\
\Omega Y_{1} \ar[r]
 & \ast \ar[r]
 & \Sigma \Omega Y_1,
}$$
which implies the condition (1).
The condition (2) is obtained as follows: for any $t\wedge x \in \Sigma L_1$, we have 
\begin{align*}
e_1^{Y_1}{\circ}\hat{f}_1(t\wedge x)
 & =
 \mathrm{ev}{\circ}\Sigma \ad{f|_{X_1}}(t\wedge x)
\\& =
 \mathrm{ev} (t\wedge \Sigma \ad{f|_{X_1}}(x))
 = \ad{f|_{X_1}}(x)(t) 
 = (f|_{X_1})(t\wedge x).
\end{align*}
\par
When $i=k$, suppose we have already obtained $\hat{f}_{i}$ and $\hat{f}^{0}_{i}$ up to $i=k{-}1$, which satisfies the conditions (1) and (2).

Firstly, we construct $\hat{f}^{0}_k : L_k\to E^k\Omega Y_k$ as follows:
the homotopy class of a map $P^{k-1}\Omega i_{k-1,k}^{Y}{\circ}\hat{f}_{k-1}{\circ}h_{k} : L_{k} \to P^{k-1}{\Omega}Y_{k}$ can be described as ${h_{k}}_{*}(P^{k-1}\Omega i_{k-1,k}^{Y}{\circ}\hat{f}_{k-1})$ $\in$ $[L_k, Y_k]$ with $P^{k-1}\Omega i_{k-1,k}^{Y}{\circ}\hat{f}_{k-1}$ $\in$ $[X_{k-1}, Y_k]$ in the following ladder of exact sequences induced from a fibre sequence $E^k\Omega Y_k \rightarrow P^{k-1}\Omega Y_k \hookrightarrow  P^{\infty}\Omega Y_{k} \ (\simeq Y_{k})$:
$$
\begin{diagram}
\node{[X_{k-1}, E^k\Omega Y_k]}
 \arrow{s,l}{h_{k}^{*}}
 \arrow{e,t}{{p_{k}^{\Omega Y_k}}_{\ast}}
\node{[X_{k-1}, P^{k-1}\Omega Y_k]}
 \arrow{s,l}{h_{k}^{*}}
 \arrow{e,t}{{e^{Y_k}_{k-1}}_{\ast}}
\node{[X_{k-1}, Y_k]}
 \arrow{s,l}{h_{k}^{*}}
\\
\node{[L_k, E^k\Omega Y_k]}
 \arrow{e,t}{{p_{k}^{\Omega Y_k}}_{\ast}}
\node{[L_k, P^{k-1}\Omega Y_k]}
 \arrow{e,t}{{e^{Y_k}_{k-1}}_{\ast}}
\node{[L_k, Y_k].}
\end{diagram}
$$
Since we know that the naturality of $e^{Z}_{k-1}$ at $Z$ implies $e_{k-1}^{Y_k}{\circ}P^{k-1}\Omega i_{k-1,k}^{Y} = i_{k-1,k}^{Y}{\circ}e_{k-1}^{Y_{k-1}}$, that the induction hypothesis implies $e_{k-1}^{Y_{k-1}}{\circ}\hat{f}_{k-1} = f\vert_{X_{k-1}}$ and that the naturally of $i^{Z}_{k-1,k}$ at $Z$ implies $i_{k-1,k}^{Y}{\circ}f|_{X_{k-1}} = f|_{X_{k}}{\circ}i_{k-1,k}^{X}$, we obtain ${e_{k-1}^{Y_k}}_{*}(P^{k-1}\Omega i_{k-1,k}^{Y}{\circ}\hat{f}_{k-1}) = i_{k-1,k}^{Y}{\circ}e_{k-1}^{Y_{k-1}}{\circ}\hat{f}_{k-1} = f|_{X_{k}}{\circ}i_{k-1,k}^{X} \in [X_{k-1},Y_{k}]$.
On the other hand, since $L_{k}\xrightarrow{h_{k-1}} X_{k-1}\hookxrightarrow{i_{k-1,k}^X} X_{k}$ is a cofibration, we obtain
$$
{e_{k-1}^{Y_k}}_{*}(h_{k}^{*}(P^{k-1}\Omega i_{k-1,k}^{Y}{\circ}\hat{f}_{k-1})) = 
f|_{X_{k}}{\circ}i_{k-1,k}^X{\circ}h_{k}=0.
$$
Thus we have 
${e_{k-1}^{Y_k}}_{\ast}(P^{k-1}\Omega i_{k-1,k}^{Y}{\circ}f_{k-1}{\circ}h_{k-1})=0$ and there exists a map $\hat{f}^{0}_k : L_k \to E^k\Omega Y_k$ such that ${p^{\Omega Y_k}_{k}}_{\ast}(\hat{f}^{0}_k)=
P^{k-1}\Omega i_{k-1,k}^Y{\circ}\hat{f}_{k-1}{\circ}h_{k-1}$.

Secondly, we define a map $f^{\prime}_k : X_k \to P^k\Omega Y_k$ as a close approximation of $\hat{f}_{k}$, which is given by
\[f'_k= P^{k-1}\Omega i_{k-1,k}^Y{\circ}\hat{f}_{k-1} \cup C(\hat{f}^{0}_k). \]
\par
Then we can easily see that $f^{\prime}_{k}$ makes the right hand square of the following diagram commutative.
$$\xymatrix{
L_k \ar[r]^{h_{k}} \ar[dd]_{\hat{f}^{0}_k}
 & X_{k-1} \ar@{^{(}->}[r]^{i_{k-1,k}^X} \ar[d]^{\hat{f}_{k-1}}
 &  X_k \ar@{-->}[dd]^{f^{\prime}_k}\\
 & P^{k-1}\Omega Y_{k-1} \ar@{^{(}->}[d]^{P^{k-1}\Omega i_{k-1,k}^Y}
 & \\
E^k\Omega Y_k \ar[r]_{p_{k}^{\Omega Y_k}}
 & P^{k-1} \Omega Y_k \ar@{^{(}->}[r]_{\iota_{k-1,k}^{\Omega Y_k}}
 & P^k \Omega Y_k. 
}$$
By its definition, $f^{\prime}_k$ satisfies the equation 
\begin{equation}\label{f'nu} 
(f^{\prime}_k\vee \Sigma \hat{f}^{0}_k){\circ}\nu_{k} = \bar{\nu}_{k}{\circ}f^{\prime}_k,
\end{equation}
where $\nu_{k} : X_k\to X_k \vee \Sigma L_k$ and
$\bar{\nu}_{k} : P^k\Omega Y_k\to P^k\Omega Y_k \vee \Sigma E^k\Omega Y_k$
are the canonical co-pairings (see Hilton \cite{Hilton65} or Oda \cite{Oda92}).
In the exact sequence 
$[X_{k-1},Y_k]\xleftarrow{{i_{k-1,k}^X}^{\ast}} [X_{k},Y_k]
  \xleftarrow{q^{\ast}}[\Sigma L_k, Y_k]$,
we have an equation
\begin{align*}&
{i_{k-1,k}^X}^{\ast}(e_k^{Y_k}{\circ}f^{\prime}_k)
  = e_k^{Y_k}{\circ}f^{\prime}_k{\circ}i_{k-1,k}^X 
  = e_k^{Y_k}{\circ}(\iota_{k-1,k}^{\Omega Y_k}{\circ}P^{k-1}\Omega i_{k-1,k}^Y
  {\circ}f_{k-1})
\\&\quad
  = e_{k-1}^{Y_k}{\circ}P^{k-1}\Omega i_{k-1,k}^Y{\circ}f_{k-1}
  = i_{k-1,k}^Y{\circ}f|_{X_k-1}
  = f|_{X_k}{\circ}i_{k-1,k}^X
  = {i_{k-1,k}^X}^{\ast}(f|_{X_k}).
\end{align*}
Thus by a standard argument of homotopy theory (see \cite{Hilton65} for example), we obtain the difference map
$\delta^{\prime}_k : \Sigma L_k \to Y_k$ such that
\[
f|_{X_k} = \nabla_{Y_k}\circ(e_k^{Y_k}{\circ}f^{\prime}_k \vee \delta^{\prime}_k){\circ}\nu_{k}.
\]
Finally we construct a map $\hat{f}_k : X_k\to P^k\Omega Y_k$.
We can observe that $\delta'_{k}$ lies in lower middle group $[\Sigma L_k, Y_k]$ in the exact sequence 
$$
\xymatrix{
\ar[r]  & [L_k, \Omega P^{k-1}\Omega Y_k]
 \ar[r]^{\hspace{15pt}{\Omega e_{k-1}^{Y_k}}_{\ast}}
 & [L_k, \Omega Y_k ]  \ar[r]^{\Delta_{\ast}\hspace{10pt}}
 & [L_k, E^k\Omega Y_k]  \ar[r] & \hspace{15pt} \\
 & [\Sigma L_k, P^{k-1}\Omega Y_k]
  \ar[u]_{\cong}^{\mathrm{ad}} \ar[r]^{\hspace{15pt}{e_{k-1}^{Y_{k}}}_{\ast}}
 & [\Sigma L_k, Y_k] \ar[u]_{\cong}^{\mathrm{ad}},
}$$
where $\Delta_{*}$ denotes the connecting map.
Since $\Omega e_{k-1}^{Y_k}$ has a section, there is a map
$\delta_k : \Sigma L_k\to P^{k-1}\Omega Y_k$ satisfying $\delta^{\prime}_k = e_{k-1}^{Y_k}{\circ}\delta_k$, and hence we have
\begin{align*}
f|_{X_k}
 & = \nabla_{Y_k}\circ(e_k^{Y_k}{\circ}f^{\prime}_k \vee e_{k-1}^{Y_k}{\circ}\delta_k)
   {\circ}\nu_{k} \\
 & = \nabla_{Y_k}\circ(e_k^{Y_k}{\circ}f^{\prime}_k \vee
 e_k^{Y_k}{\circ}\iota_{k-1,k}^{\Omega Y_k}{\circ}\delta_k){\circ}\nu_{k} \\
 & = \nabla_{Y_k}\circ(e_k^{Y_k}\vee e_k^{Y_k}){\circ}(f^{\prime}_k \vee \iota_{k-1,k}^{\Omega Y_k}\circ\delta_k)
{\circ}\nu_{k} \\
 & = e_k^{Y_k}{\circ}\nabla_{P^k\Omega{Y_k}}{\circ}(f^{\prime}_k \vee \iota_{k-1,k}^{\Omega Y_k}\circ\delta_k)
{\circ}\nu_{k}. \\
\end{align*}
So we define the map $\hat{f}_k$ to be $\nabla_{P^k\Omega{Y_k}}{\circ}(f^{\prime}_k \vee \iota_{k-1,k}^{\Omega Y_k}\circ\delta_k){\circ}\nu_{k}$.
Then we can easily check that the condition (2) for $i=k$.
Further, since $\nu_{k}$ is a co-pairing, we have
$$pr_1{\circ}\nu_{k}{\circ}i_{k-1,k}^X= \mathrm{id}_{X_{k}}{\circ}i_{k-1,k}^X = i_{k-1,k}^X
\text{\ \ and\ \ }
pr_2{\circ}\nu_{k}{\circ}i_{k-1,k}^X= q{\circ}i_{k-1,k}^X = 0,$$
where $pr_1 : X_k \vee \Sigma L_k \to X_k$ and 
$pr_2 : X_k \vee \Sigma L_k \to \Sigma L_k$ are the first and second projections,
respectively.
Hence we obtain the equation
\begin{align*}
\hat{f}_k{\circ}i_{k-1,k}^X
 &
  = \nabla_{P^k\Omega{Y_k}}{\circ}
  (f^{\prime}_k \vee \iota_{k-1,k}^{\Omega Y_k}\circ\delta_k){\circ}\nu_{k}{\circ}i_{k-1,k}^X
 \\&
  = f^{\prime}_k{\circ}i_{k-1,k}^X
  = \iota_{k-1,k}^{\Omega Y_k}{\circ}P^{k-1}\Omega i_{k-1,k}^Y{\circ}\hat{f}_{k-1},
\end{align*}
which implies the condition (1) for $i=k$.
This completes the induction step and we obtain the proposition.
\end{proof}

Let $\{\hat{f}_{i} : X_{i}\to P^{i}\Omega Y_{i} \,|\, 0{\leq} i {\leq} n \}$ and 
 $\{\hat{f}^{0}_{i} : L_{i}\to E^{i}\Omega Y_{i} \,|\,1{\leq} i {\leq} m \}$ be maps obtained from the filtered map $f : X \to Y$ as in Proposition \ref{cd-to-gs}.
Let $\nu_{i} : X_i\to X_i \vee \Sigma L_i$ and
$\bar{\nu}_{i} : P^i\Omega Y_i\to P^k\Omega Y_i \vee \Sigma E^i\Omega Y_i$ be the canonical co-pairings.
\begin{prop}\label{coop-commu}
If the complex $L_i$ be a co-H-space, then the following diagram is commutative. 
$$\xymatrix{
X_i \ar[r]^{\nu_{i} } \ar[d]_{\hat{f}_i}
 & X_k \vee \Sigma L_i  \ar[d]^{\hat{f}_i \vee \Sigma \hat{f}^{0}_i} \\
P^i \Omega Y_i \ar[r]_{\bar{\nu}_{i} \hspace{25pt}}
 & P^i \Omega Y_i \vee \Sigma E^i \Omega Y_i.
}$$
\end{prop}
\begin{proof}
By the definition of $\hat{f}_i$ using $f'_{i}$ in the proof of Proposition \ref{cd-to-gs}, and by the relation between the composition and the wedge of maps, we have
\begin{align*}
(\hat{f}_i \vee \Sigma \hat{f}^{0}_i){\circ}\nu_{i}
 & = \{(\nabla_{P}{\circ}(f^{\prime}_i \vee \iota_{i-1,i}^{\Omega Y_i}\circ\delta_i){\circ}\nu_{i})
       \vee \Sigma \hat{f}^{0}_i \}{\circ}\nu_{i} \\
 & = \{\nabla_{P}{\circ}(f^{\prime}_i \vee \iota_{i-1,i}^{\Omega Y_i}{\circ}\delta_i)
       \vee \Sigma \hat{f}^{0}_i \}{\circ}(\nu_{i} \vee \mathrm{id}_{\Sigma L_i}){\circ}\nu_{i} \\
 & = (\nabla_{P} \vee \mathrm{id}_{E})
     {\circ}(f^{\prime}_i\vee \iota_{i-1,i}^{\Omega Y_i}{\circ}\delta_i \vee \Sigma \hat{f}^{0}_i )
     {\circ}(\nu_{i} \vee \mathrm{id}_{\Sigma L_i}){\circ}\nu_{i},
\end{align*}
where $\nabla_{P}=\nabla_{P^i\Omega{Y_i}}$ and $\mathrm{id}_{E}=\mathrm{id}_{\Sigma E^i\Omega Y_i}$.
Since $L_i$ is the co-H-space, we have the equations
$$\upsilon_{i}= T\circ\upsilon_{i}
\text{\ \ and \ }
(\nu_{i} \vee \mathrm{id}_{\Sigma L_i}){\circ}\nu_{i}
   = (\mathrm{id}_{X_i} \vee \upsilon_{i}){\circ}\nu_{i},
$$
where $\upsilon_{i} : \Sigma L_i \to \Sigma L_i \vee \Sigma L_i$ is the co-multiplication and
$T : \Sigma L_i \vee \Sigma L_i \to \Sigma L_i \vee \Sigma L_i$
is the commutative map.
So we can proceed as follows:
\begin{align*}
(\hat{f}_{i} \vee \Sigma \hat{f}^{0}_{i}){\circ}\nu_{i}
 & = (\nabla_{P} \vee \mathrm{id}_{E})
     {\circ}(f^{\prime}_{i}\vee \iota_{i-1,i}^{\Omega Y_{i}}{\circ}\delta_{i} \vee \Sigma \hat{f}^{0}_{i} )
     {\circ}(\mathrm{id}_{X_{i}} \vee \upsilon_{i}){\circ}\nu_{i} \\
 & = (\nabla_{P} \vee \mathrm{id}_{E})
     {\circ}(f^{\prime}_{i}\vee \iota_{i-1,i}^{\Omega Y_{i}}{\circ}\delta_{i} \vee \Sigma \hat{f}^{0}_{i} )
     {\circ}(\mathrm{id}_{X_{i}} \vee T\circ\upsilon_{i}){\circ}\nu_{i} \\
 & = (\nabla_{P} \vee \mathrm{id}_{E})
     {\circ}\{f^{\prime}_{i} \vee T^{\prime}{\circ}
     (\Sigma \hat{f}^{0}_{i} \vee \iota_{i-1,i}^{\Omega Y_{i}}{\circ}\delta_{i})\}
     {\circ}(\mathrm{id}_{X_{i}} \vee \upsilon_{i} ){\circ}\nu_{i} \\
 & = (\nabla_{P} \vee \mathrm{id}_{E})
     {\circ}(f^{\prime}_{i}\vee T^{\prime} )\\
 &\hspace{50pt}{\circ}\{\mathrm{id}_{X_{i}}\vee (\Sigma \hat{f}^{0}_{i} \vee \iota_{i-1,i}^{\Omega Y_{i}}{\circ}\delta_{i})\}
     {\circ}(\mathrm{id}_{X_{i}}\vee \upsilon_{i} ){\circ}\nu_{i}\\
 & = (\nabla_{P} \vee \mathrm{id}_{E})
     {\circ}(\mathrm{id}_{P} \vee T^{\prime})\\
 &\hspace{50pt}{\circ}(f^{\prime}_{i} \vee \Sigma \hat{f}^{0}_{i} \vee \iota_{i-1,i}^{\Omega Y_{i}}{\circ}\delta_{i})
     {\circ}(\nu_{i} \vee \mathrm{id}_{\Sigma L_{i}}){\circ}\nu_{i}\\
 & = (\nabla_{P} \vee \mathrm{id}_{E})
     {\circ}(\mathrm{id}_{P} \vee T^{\prime})
     {\circ}\{(f^{\prime}_{i} \vee \Sigma \hat{f}^{0}_{i}){\circ}\nu_{i} \vee \iota_{i-1,i}^{\Omega Y_{i}}{\circ}\delta_{i} \}{\circ}\nu_{i},
\end{align*}
where $T^{\prime} : \Sigma E^{i}\Omega Y_{i}\vee P^{i}\Omega{Y_{i}}
    \to P^{i}\Omega{Y_{i}}\vee\Sigma E^{i}\Omega Y_{i}$
is the commutative map and $\mathrm{id}_{P}= \mathrm{id}_{P^{i}\Omega{Y_{i}}}$.
By the equation (\ref{f'nu}), we proceed further as follows:
\begin{align*}
(\hat{f}_{i} \vee \Sigma \hat{f}^{0}_{i}){\circ}\nu_{i}
 &
  = (\nabla_{P} \vee \mathrm{id}_{E})
     {\circ}(\mathrm{id}_{P} \vee T^{\prime})
     {\circ}\{(\bar{\nu}_{i}{\circ}f^{\prime}_{i} \vee \iota_{i-1,i}^{\Omega Y_{i}}{\circ}\delta_{i} \}
     {\circ}\nu_{i}
 \\&
  = (\nabla_{P} \vee \mathrm{id}_{E})
    {\circ}(\mathrm{id}_{P} \vee T^{\prime})
    {\circ}(\bar{\nu}_{i}  \vee \iota_{i-1,i}^{\Omega Y_{i}})
    {\circ}(f^{\prime}_{i} \vee \delta_{i})
    {\circ}\nu_{i}
 \\&
  = (\nabla_{P} \vee \nabla_{\Sigma E^i\Omega Y_i})
    {\circ}(\mathrm{id}_{P} \vee T^{\prime} \vee \mathrm{id}_{E})
 \\&
  \hspace{50pt}{\circ}(\bar{\nu}_{i}  \vee \bar{\nu}_{i})
    {\circ}(\mathrm{id}_{P} \vee\iota_{i-1,i}^{\Omega Y_i})
    {\circ}(f^{\prime}_i \vee \delta_i)
    {\circ}\nu_{i}
 \\&
  = \bar{\nu}_{i}{\circ}\nabla_{P}
    {\circ}(f^{\prime}_i\vee \iota_{i-1,i}^{\Omega Y_i}{\circ}\delta_i )
    {\circ}\nu_{i} 
  = \bar{\nu}_{i}{\circ}f_i.
\end{align*}
This implies the proposition.
\end{proof}
\section{Cone-Decomposition associated with projective spaces}%
\label{cone-decom}
As usual, $X^{(k)}$ stands for the $k$-skeleton of a space $X$ and $f^{(k)} : X^{(k)} \to Y^{(k)}$ stands for $f\vert_{X^{(k)}}$, since $f\vert_{X^{(k)}} : X^{(k)} \to Y$ is compressible into $Y^{(k)}$.
When the dimension of $X$ is less than or equal to $n$, 
then we do not distinguish $f^{(n)} : X \to Y^{(n)}$ and $f : X\to Y$, by the same reason.
\par
Let $G$ be a compact Lie group with a cone-decomposition of length $m$,
that is, there is a series of cofibre sequences
\begin{equation}
\label{cd-fm}
\{K_i\xrightarrow{h_{i}} F_{i-1}\hookrightarrow F_i \,|\,1 \leq i \leq m \}
\end{equation}
with $F_0=\{\ast\}$ and $F_m \simeq G$.
We also name the inclusion as $i^{F}_{i-1,i} : F_{i-1}\hookrightarrow F_i$ and its quotient as $q^{F}_{i-1,i} : F_{i} \to \Sigma{K_{i}}$.
Let $\ell$ be the dimension of Lie group $G$.
\begin{lem}\label{cofibPmFm}
Suppose that the complex $K_m$ is the sphere $S^{\ell-1}$ and $\ell\geq 3$, $m\geq 3$.
Then there is a cofibre sequence as follows:
\[
(E^{m}\Omega F_{m-1})^{(\ell-1)}\vee K_m\xrightarrow{p^{\prime}}
 (P^{m-1}\Omega F_{m-1})^{(\ell)}\hookrightarrow (P^m\Omega F_m)^{(\ell)}.
\]
\end{lem}
\begin{proof}
First, we determine the homotopy type of the $(\ell{-}1)$-skeleton of the homotopy fibre of the map
$P^{m-1}\Omega i^{F}_{m-1,m} : P^{m-1}\Omega F_{m-1}\to P^{m-1}\Omega F_m$.
Let $\mathfrak{F}$ be the homotopy fibre of $P^{m-1}\Omega i^{F}_{m-1,m}$, which fits in with the following commutative diagram whose rows and columns are all fibrations.
$$\xymatrix{
   \Omega (E^m\Omega F_m, E^m\Omega F_{m-1})
     \ar[d]
     \ar[r]
 & \mathfrak{F}
     \ar[d]
     \ar[r]
 & \Omega (F_m, F_{m-1})
     \ar[d]
\\
   E^m\Omega F_{m-1}
     \ar@{^{(}->}[d]^{E^m\Omega i^{F}_{m-1,m}}
     \ar[r]^{p_{m}^{\Omega F_{m-1}}}
 & P^{m-1}\Omega F_{m-1}
     \ar@{^{(}->}[d]^{P^{m-1}\Omega i^{F}_{m-1,m}}
     \ar[r]^{e^{F_{m-1}}_{m-1}}
 & F_{m-1}
     \ar@{^{(}->}[d]^{i^{F}_{m-1,m}}
\\
   E^m\Omega F_m
     \ar[r]^{p_{m}^{\Omega F_{m}}}
 & P^{m-1}\Omega F_m
     \ar[r]^{e^{F_{m}}_{m-1}}
 & F_m.
}$$
Since the pair $(F_m, F_{m-1})$\hspace{-1.5pt} is $(\ell{-}1)$-connected, 
$(\Omega F_m, \Omega F_{m-1})$
 is $(\ell{-}2)$-connected 
and $(E^m\Omega F_m, E^m\Omega F_{m-1})$ is $(\ell{+}m{-}3)$-connected.
Then $\Omega (E^m\Omega F_m, E^m\Omega F_{m-1})$ is $(\ell{+}m{-}4)$-connected and hence $(\ell{-}1)$-connected.
By using homotopy exact sequence, we obtain that $\mathfrak{F}$ is $(\ell{-}2)$-connected.
Thus $\mathfrak{F}$ is $1$-connected.

On the other hand, by using Serre exact sequence 
\[
H_{2\ell+m-5}(\Omega (E^m\Omega F_m, E^m\Omega F_{m-1}))\to 
\cdots \to H_{k}(\Omega (E^m\Omega F_m, E^m\Omega F_{m-1}))
  \to 
\]
\[
H_{k}(\mathfrak{F}) 
\to H_{k}(\Omega (F_m, F_{m-1}))\to H_{k-1}(\Omega (E^m\Omega F_m, E^m\Omega F_{m-1}))\to\cdots ,
\]
we obtain that $H_{k}(\mathfrak{F})$ is isomorphic to $H_{k}(\Omega (F_m, F_{m-1}))$ for $k \leq \ell{+}m{-}3 \ (\geq \ell)$.
By Blakers-Massey's theorem, we obtain $\pi_l(F_m, F_{m-1}) \cong \pi_l (S^l)$, and hence 
\[
\pi_{\ell-1} (\Omega(F_m, F_{m-1})) \cong \pi_l (F_m, F_{m-1}) \cong \pi_l (S^l) \cong \mathbb{Z}.
\]
Then by Hurewicz Isomorphism Theorem, we obtain
\[
H_{\ell-1}(\mathfrak{F}) \cong H_{\ell-1}(\Omega (F_m, F_{m-1}))\cong \pi_{\ell-1} (\Omega(F_m, F_{m-1}))\cong \mathbb{Z}.
\]
Thus $\mathfrak{F}$ has the following homology decomposition.
$$
\mathfrak{F} \simeq S^{\ell-1} \cup \ \text{(Moore cells in dimensions $\geq \ell$)}.
$$
By Ganea's fibre-cofibre construction (see Ganea \cite{Ganea67}), we obtain a map 
$$
\phi_{0} : P^{m-1}\Omega F_{m-1} \cup C\mathfrak{F} \to P^{m-1}\Omega F_{m},
$$
as the homotopy pushout 
$$
\begin{diagram}
\node{\mathfrak{F}}
     \arrow[2]{s}
     \arrow[2]{e}
\node{}
\node{P^{m-1}\Omega F_{m-1}}
     \arrow[2]{s,J}
\\
\node{}
\node{\mbox{HPO}}
\\
\node{\{\ast\}}
     \arrow[2]{e}
\node{}
\node{P^{m-1}\Omega F_{m-1} \cup C\mathfrak{F},}
\end{diagram}
$$
which has the homotopy type of the homotopy pullback of the diagonal 
$$
\Delta : P^{m-1}\Omega F_{m} \to P^{m-1}\Omega F_{m}{\times}P^{m-1}\Omega F_{m}
$$
and the inclusion 
$$
P^{m-1}\Omega F_{m-1}{\times}P^{m-1}\Omega F_{m}{\cup}P^{m-1}\Omega F_{m}{\times}\{\ast\} \hookrightarrow P^{m-1}\Omega F_{m}{\times}P^{m-1}\Omega F_{m}:
$$
$$
\begin{diagram}
\node{P^{m-1}\Omega F_{m-1} \cup C\mathfrak{F}}
     \arrow[2]{s,l}{\phi_{0}}
     \arrow[2]{e}
\node{}
\node{\mbox{$\begin{array}{c}P^{m-1}\Omega F_{m-1}{\times}P^{m-1}\Omega F_{m} \\ \cup P^{m-1}\Omega F_{m}{\times}\{\ast\}\end{array}$}}
     \arrow[2]{s,J}
\\
\node{}
\node{\mbox{HPB}}
\\
\node{P^{m-1}\Omega F_m}
     \arrow[2]{e,b}{\Delta}
\node{}
\node{P^{m-1}\Omega F_{m}{\times}P^{m-1}\Omega F_{m}.}
\end{diagram}
$$
(see, for example, \cite[Lemma 2.1]{Iwase98} with $(X,A)=(P^{m-1}\Omega F_{m},P^{m-1}\Omega F_{m-1})$, $(Y,B)=(P^{m-1}\Omega F_{m},\{\ast\})$ and $Z=P^{m-1}\Omega F_{m}$).
Hence $\mathfrak{F}_{0}$ is given by the pullback of the trivial map
$$
\{\ast\} \to P^{m-1}\Omega F_{m}{\times}P^{m-1}\Omega F_{m}
$$
and the inclusion 
$$
P^{m-1}\Omega F_{m-1}{\times}P^{m-1}\Omega F_{m}{\cup}P^{m-1}\Omega F_{m}{\times}\{\ast\} \hookrightarrow P^{m-1}\Omega F_{m}{\times}P^{m-1}\Omega F_{m}
$$
which has the homotopy type of the pushout 
$$
\begin{diagram}
\node{\mathfrak{F}{\times}\Omega P^{m-1}\Omega F_{m}}
     \arrow[2]{s}
     \arrow[2]{e}
\node{}
\node{P^{m-1}\Omega F_{m-1}}
     \arrow[2]{s}
\\
\node{}
\node{\mbox{HPO}}
\\
\node{\mathfrak{F}}
     \arrow[2]{e}
\node{}
\node{\mathfrak{F}_{0}.}
\end{diagram}
$$
(see, for example, \cite[Lemma 2.1]{Iwase98} with $(X,A)=(P^{m-1}\Omega F_{m},P^{m-1}\Omega F_{m-1})$, $(Y,B)=(P^{m-1}\Omega F_{m},\{\ast\})$ and $Z=\{\ast\}$).
Thus the homotopy fibre $\mathfrak{F}_{0}$ of $\phi_{0}$ has the homotopy type of the join $\mathfrak{F}{\ast}\Omega P^{m-1}\Omega F_{m}$ and is ($\ell{-}1$)-connected, and hence $\phi_{0}$ is $\ell$-connected.
Thus we have that 
$$
(P^{m-1}\Omega F_{m-1})^{(\ell)} \cup CS^{\ell-1} \simeq (P^{m-1}\Omega F_{m})^{(\ell)}.
$$
We are now ready to show that $(P^m\Omega F_{m-1})^{(\ell)}\cup CS^{\ell-1} \simeq (P^m\Omega F_m)^{(\ell)}$:
since $(E^m\Omega F_m, E^m\Omega F_{m-1})$ is $(\ell{+}m{-}3)$-connected and $m {\geq} 3$, we have $(E^m\Omega F_{m-1})^{(\ell-1)}$ $\simeq$ $(E^m\Omega F_m)^{(\ell-1)}$ and hence we obtain
\begin{align*}
(P^{m}\Omega F_{m-1})^{(\ell)} \cup CS^{\ell-1} 
&\simeq 
(P^{m-1}\Omega F_{m-1})^{(\ell)} \cup C(S^{\ell-1} \vee (E^m\Omega F_{m-1})^{(\ell-1)})
\\&\simeq 
(P^{m-1}\Omega F_{m})^{(\ell)} \cup C(E^m\Omega F_{m-1})^{(\ell-1)} 
\simeq
(P^{m}\Omega F_{m})^{(\ell)}.
\end{align*}
This completes the proof of Lemma \ref{cofibPmFm}.
\end{proof}
Using Lemma \ref{cofibPmFm}, we construct cone-decompositions of $F_m{\times} F'_1$, %
$(P^m\Omega F_m)^{(\ell)}$ and $(P^m\Omega F_m)^{(\ell)}{\times}(\Sigma \Omega F'_1)^{(\ell)}$:
first, we construct a cone-decomposition of $F_m{\times}F'_1$. 
Let $K^{m,1}_i$ and $F^{m,1}_i$ be as follows.
\begin{align*}&
K^{m,1}_i = K_i \vee \{K_{i-1}{\ast} K'_1\}
& \text{for \ } 1\leq i\leq m,
\\&
F^{m,1}_i = F_i{\times} \{\ast\}\cup F_{i-1}{\times} F'_1
& \text{for \ } 0\leq i\leq m,
\\&
K^{m,1}_{m+1} = K_m{\ast}K'_1 \text{\ \ and\ \ }F^{m,1}_{m+1}=F_m{\times} F'_1,
\end{align*}
where $K_{0}$ and $F_{-1}$ are empty sets.
We define ${w}^{m,1}_{i} : K^{m,1}_{i} \to F^{m,1}_{i-1}$ by
\begin{align*}&
{w}^{m,1}_{i}\vert_{K_{i}} = (\incl){\circ}(h_{i}{\times} \{\ast\}) : K_i \to F_{i-1} = F_{i-1}{\times} \{\ast\}
\hookrightarrow F^{m,1}_{i-1}
\quad\text{and}
\\&
{w}^{m,1}_{i}\vert_{K_{i-1}{\ast} K'_1} = [\chi_{i-1}, \Sigma\mathrm{id}_{K'_1}]^r
  : K_{i-1}{\ast} K'_1 \to F_{i-1}{\times} \{\ast\}\cup F_{i-2}{\times} \Sigma K'_1 = F^{m,1}_{i-1}
\intertext{for $1\leq i\leq m$, and}&
{w}^{m,1}_{m+1} = [\chi_m , \Sigma\mathrm{id}_{K'_1}]^r : K^{m,1}_{m+1} \to F^{m,1}_{m},
\end{align*}
where $[\chi_i , \Sigma\mathrm{id}_{K'_1}]^r$ denotes the relative Whitehead product of the characteristic map $\chi_i : (CK_{i},K_{i})\to (F_{i}, F_{i-1})$ and the suspension of the identity map $\Sigma\mathrm{id}_{K'_1}$.
Let
$i^{m,1}_{i} : F^{m,1}_{i}\to F^{m,1}_{i+1}$ be the canonical inclusion 
 for $0\leq i\leq m$.
Then the sequence of cofibrations
\begin{equation}
\label{cd-FxF}
\{ K^{m,1}_{i} \xrightarrow{{w}^{m,1}_{i}}
    F^{m,1}_{i-1} \hookxrightarrow{i^{m,1}_{i-1}} F^{m,1}_{i}
 \,|\, 1\le i \le m+1 \}
\end{equation}
is a cone-decomposition of 
$F_m {\times} F'_1$
of length $m+1$.
\par
Second, we construct a cone-decomposition of
$(P^m\Omega F_m)^{(\ell)}$.
By lemma \ref{cofibPmFm}, we obtain a cone-decomposition of 
$(P^m\Omega F_m)^{(\ell)}$ of length $m$.
$$
\left\{
\begin{array}{l}
(\Omega F_{m-1})^{(\ell-1)}\to \{\ast \}
  \hookrightarrow (\Sigma \Omega F_{m-1})^{(\ell)}
\\[2mm]
(E^2\Omega F_{m-1})^{(\ell-1)}\to (\Sigma \Omega F_{m-1})^{(\ell)}
  \hookrightarrow (P^2\Omega F_{m-1})^{(\ell)}
\\[2mm]
\hspace{100pt} \vdots
\\[2mm]
(E^{m-1}\Omega F_{m-1})^{(\ell-1)} \to (P^{m-2}\Omega F_{m-1})^{(\ell)}
  \hookrightarrow (P^{m-1}\Omega F_{m-1})^{(\ell)}
\\[2mm]
(E^m\Omega F_{m-1})^{(\ell-1)}\vee K_m \to (P^{m-1}\Omega F_{m-1})^{(\ell)}
  \hookrightarrow (P^m\Omega F_m)^{(\ell)}.
\end{array}
\right.
$$
\par
Third, we construct a cone-decomposition of
$(P^m\Omega F_m)^{(\ell)}{\times} (\Sigma \Omega F'_1)^{(\ell)}$.
We define two series of spaces $\hat{E}_i$ and $\hat{F}_i$, $1 \leq i \leq m{+}1$ as follows.
\begin{align*}&
\hat{E}_{i}
= (E^i\Omega F_{m-1})^{(\ell-1)} 
   \vee\{(E^{i-1}\Omega F_{m-1})^{(\ell-1)}{\ast} {(\Omega F'_1)}^{(\ell-1)}\},
\\&
\hat{F}_{i}\hspace{5pt}
= (P^{i}\Omega F_{m-1})^{(\ell)} {\times} \{\ast\}
   \cup (P^{i-1}\Omega F_{m-1})^{(\ell)} {\times} (\Sigma\Omega F'_1)^{(\ell)}
\intertext{for $1\leq i < m$,}&
\hat{E}_{m}
= \{(E^m\Omega F_{m-1})^{(\ell-1)}\vee K_m\} 
   \vee\{(E^{m-1}\Omega F_{m-1})^{(\ell-1)}{\ast} {(\Omega F'_1)}^{(\ell-1)}\},
\\&
\hat{F}_{m}
= (P^{m}\Omega F_{m})^{(\ell)} {\times} \{\ast\}
   \cup (P^{m-1}\Omega F_{m-1})^{(\ell)} {\times} (\Sigma\Omega F'_1)^{(\ell)}
\intertext{for $i=m$ and, for $i=m{+}1$,}&
\hat{E}_{m+1}
= \{(E^{m}\Omega F_{m-1})^{(\ell-1)} \vee K_m\}{\ast} (\Omega F'_1)^{(\ell-1)},
\\&
\hat{F}_{m+1}
= (P^{m}\Omega F_{m})^{(\ell)} {\times} (\Sigma\Omega F'_1)^{(\ell)}.
\end{align*}
\begin{rem}
$E^{-1}\Omega F_{m-1}$ and $P^{-1}\Omega F_{m-1}$
denote empty sets.
\end{rem}
We further define maps $\hat{w}_{i} : \hat{E}_{i} \to \hat{F}_{i-1}$, $1 \leq i \leq m{+}1$ by
\begin{align*}&
\hat{w}_{i}\vert_{(E^i\Omega F_{m-1})^{(\ell-1)}} = (\incl){\circ}(p_{i}^{\Omega F_{m-1}})^{(\ell-1)}
\\&\qquad\qquad\qquad
	: (E^i\Omega F_{m-1})^{(\ell-1)} 
    \to (P^{i-1}\Omega F_{m-1})^{(\ell)} \hookrightarrow \hat{F}_{i-1}
\\&
\hat{w}_{i}\vert_{(E^{i-1}\Omega F_{m-1})^{(\ell-1)}{\ast} {(\Omega F'_1)}^{(\ell-1)}} = [\chi^{\prime}_{i-1} , (\mathrm{id}_{\Sigma\Omega F'_1})^{(\ell-1)}]^r
\\&\qquad\qquad\qquad
	: (E^{i-1}\Omega F_{m-1})^{(\ell-1)}{\ast} {(\Omega F'_1)}^{(\ell-1)}
   \to \hat{F}_{i-1}
\intertext{for $1\leq i< m$,}&
\hat{w}_{m} = (\incl){\circ}p^{\prime}
	: (E^m\Omega F_{m-1})^{(\ell-1)}\vee K_m \to
   (P^{m-1}\Omega F_{m-1})^{(\ell)} \hookrightarrow \hat{F}_{m-1}
\\&
\hat{w}_{m} = [\chi^{\prime}_{m-1} , (\mathrm{id}_{\Sigma\Omega F'_1})^{(\ell-1)}]^r : (E^{m-1}\Omega F_{m-1})^{(\ell-1)}{\ast} {(\Omega F'_1)}^{(\ell-1)}
   \to \hat{F}_{m-1}, 
\intertext{for $i=m$ and, for $i=m{+}1$,}&
\hat{w}_{m+1} = [\chi^{\prime}_m , (\mathrm{id}_{\Sigma\Omega F'_1})^{(\ell-1)}]^r : \hat{E}_{m+1} \to \hat{F}_{m},
\end{align*}
where $p^{\prime} : (E^m\Omega F_{m-1})^{(\ell-1)}\vee K_m \to (P^{m-1}\Omega F_{m-1})^{(\ell)}$
is the map given in Lemma \ref{cofibPmFm} and $\chi^{\prime}_{i}$, $1 \leq i \leq m$, are the following characteristic maps.
\begin{align*}&
\chi^{\prime}_i : (C(E^{i}\Omega F_{m-1})^{(\ell-1)}, (E^{i}\Omega F_{m-1})^{(\ell-1)})
 \to ((P^{i}\Omega F_{m-1})^{(\ell)}, (P^{i-1}\Omega F_{m-1})^{(\ell)})
\intertext{for $0\leq i< m$ and for $i=m$,}&
\chi^{\prime}_{m} : (CE^{\prime},E^{\prime})
 \to ((P^{m}\Omega F_{m-1})^{(\ell)}, (P^{m-1}\Omega F_{m-1})^{(\ell)}),
\end{align*}
and $E^{\prime} = (E^{m}\Omega F_{m-1})^{(\ell-1)} \vee K_m$.
We denote $\hat{i}_{i} : \hat{F}_{i}\to\hat{F}_{i+1}$ by the canonical inclusion for $0\leq i\leq m$.
Then the sequence of cofibrations
\begin{equation}
\label{cd-pxp}
\{ \hat{E}_{i} \xrightarrow{\hat{w}_{i}}
    \hat{F}_{i-1} \hookxrightarrow{\hat{i}_{i-1}} \hat{F}_{i}
 \,|\, 1\le i \le m+1 \}
\end{equation}
gives a cone-decomposition of 
$(P^m\Omega F_m)^{(\ell)}{\times} (\Sigma \Omega F'_1)^{(\ell)}$
of length $m+1$.
\section{Structure map and cone-decomposition}
\label{str-cone}
We fix the structure maps $\sigma_{i}$ for $\cat{F_{i}} \leq i$ and $\sigma'$ for $\cat{F'_{1}} \leq 1$ as follows:
\begin{defi}
The identity map $\id{F_{m}} : F_m\to F_m$ is a filtered map w.r.t. the filtration $\ast = F_{0} \subset F_{1} \subset \cdots \subset F_{m}$.
Then by Proposition \ref{cd-to-gs}, we obtain maps $\sigma_i = \widehat{\id{F_{m}}}_{i} : F_i\to P^i\Omega F_i$ for $1\leq i\leq m$ and $\widehat{\id{F_{m}}}^{0}_j : K_j\to E^j\Omega F_j$ for $1\leq j\leq m$.
Let $g_j = \widehat{\id{F_{m}}}^{0}_j : K_j \to (E^j\Omega F_j)^{(\ell-1)}$ for $1\leq j\leq m{-}1$
and 
$g_m$ the composition $K_m \xrightarrow{\widehat{\id{F_{m}}}^{0}_m} (E^{m}\Omega F_{m})^{(\ell-1)} = 
 (E^{m}\Omega F_{m-1})^{(\ell-1)} \hookrightarrow 
 (E^{m}\Omega F_{m-1})^{(\ell-1)} \vee K_m$.
We also obtain $g'=\ad{\mathrm{id}_{K'_1}} : K'_1\to \Omega\Sigma K'_1=\Omega F'_1$ and
$\sigma'=\Sigma g' : F'_1 \to \Sigma\Omega F'_1$.
\end{defi}
\par
Let $\nu^{m,1}_{k} : F^{m,1}_{k}\to F^{m,1}_{k}\vee\Sigma K^{m,1}_{k}$
and $\hat{\nu}_{k} : \hat{F}_{k}\to \hat{F}_{k}\vee\Sigma \hat{K}_{k}$
be the canonical co-pairings for $1\le k\le m{+}1$.
\begin{lem}\label{2cofib}
The following diagram is commutative.
$$\xymatrix{
K^{m,1}_{m+1}
   \ar[r]^{{w}^{m,1}_{m+1}} \ar[d]^{g_m*g' }
 & F^{m,1}_m \ar[r]^{i^{m,1}_m}
    \ar[d]^{\hat{\sigma}_m}
 & F^{m,1}_{m+1} \ar[rr]^{\hspace{-20pt}\nu^{m,1}_{m+1}} \ar[d]^{\sigma_m{\times} \sigma'}
 &
 & F^{m,1}_{m+1} \vee \Sigma K^{m,1}_{m+1}
   \ar[d]^{\sigma_m{\times} \sigma_1 \vee \Sigma g_m{\ast} g'}
\\
\hat{E}_{m+1} \ar[r]^{\hat{w}_{m+1}}
 & \hat{F}_{m} \ar[r]^{\hat{i}_m}
 & \hat{F}_{m+1} \ar[rr]^{\hspace{-20pt}\hat{\nu}_{m+1}}
 &
 & \hat{F}_{m+1} \vee \Sigma \hat{E}_{m+1}.
}$$
Here the map $\hat{\sigma}_m = \sigma_m{\times}\{\ast\}\cup \sigma_{m-1}{\times} \sigma'$.
\end{lem}
To prove this Lemma,  we need to show the following propositions.
\begin{prop}\label{prop:key-prop}
Let $K \xrightarrow{f} A \hookrightarrow C(f)$ and $L \xrightarrow{g} B \hookrightarrow C(g)$ be two cofibre sequences with canonical co-pairings $\nu_{f} : C(f) \to C(f) \vee \Sigma{K}$ and $\nu_{g} : C(g) \to C(g) \vee \Sigma{L}$, respectively.
Then the canonical co-pairing $\nu : C(f){\times}C(g) \to C(f){\times}C(g) \vee \Sigma{K{\ast}L}$ is given by the following composition of maps.
\begin{align*}&
C(f){\times}C(g) 
\\&\quad\xrightarrow{\nu_{f}{\times}\nu_{g}}  C(f){\times}C(g) \cup_{C(f)} C(f){\times}\Sigma{L} \cup_{C(g)} \Sigma{K}{\times}C(g) \cup_{\Sigma{K}\vee\Sigma{L}} \Sigma{K}{\times}\Sigma{L} 
\\&\qquad
\xrightarrow{~\Phi~} C(f){\times}C(g) \vee \Sigma{K}{\times}\Sigma{L}/(\Sigma{K}\vee\Sigma{L}) \xrightarrow{\approx} C(f){\times}C(g) \vee \Sigma(K{\ast}L),
\end{align*}
\begin{figure}[!ht]
\includegraphics[width=12cm]{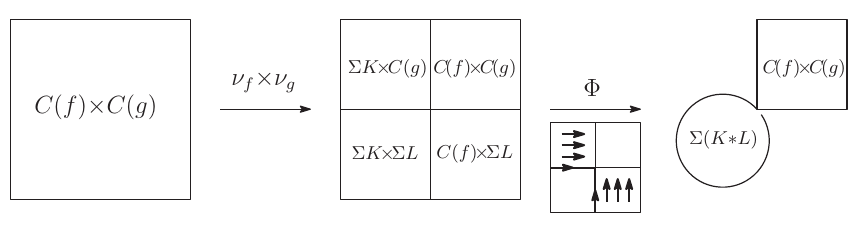}
\caption{}\label{fig-phi}
\end{figure}\par\noindent
where $\Phi$ is given by $\Phi\vert_{C(f){\times}\Sigma{L}} = \proj_{1}$ the first projection, $\Phi\vert_{\Sigma{K}{\times}C(g)} = \proj_{2}$ the second projection, and $\Phi\vert_{\Sigma{K}\times\Sigma{L}} : \Sigma{K}\times\Sigma{L} \twoheadrightarrow \Sigma{K}{\times}\Sigma{L}/(\Sigma{K}{\vee}\Sigma{L})$ the canonical collapsing.
\end{prop}
\begin{proof}
We give here explicit descriptions of $C_{h}$ for a map $h : X \to Z$ and related spaces as follows.
\begin{align*}&
CX = ([0,1]{\times}X) \amalg \{\ast\}/\sim, \quad (0,x) \sim \ast;
\\&
C_{<\frac{1}{2}}X = \{t{\wedge}x \in CX \,\vert\, t < \frac{1}{2}\} 
\subset 
C_{{\leq}\frac{1}{2}}X = \{t{\wedge}x \in CX \,\vert\, t \leq \frac{1}{2}\},
\\&
C(h) = Z \amalg CX/\sim, \quad 1{\wedge}x \sim f(x)\quad (\text{mapping cone of $h$}),
\\&
C_{\geq\frac{1}{2}}(h) = \{t{\wedge}x \in C(h) \,\vert\, t \geq \frac{1}{2}\} \quad (\text{mapping cylinder of $h$}),
\end{align*}
where we denote the class of $(t,x) \in [0,1]{\times}X$ in $CX$ or $C(h)$ by $t{\wedge}x$.
Using these notions, we describe $C(f){\times}C(g)$ as the cofibre of a Whitehead product $[\chi_{f},\chi_{g}] : K{\ast}L = CK{\times}L \cup K{\times}CL \to C(f){\times}B \cup A{\times}C(g)$ of $\chi_{f} : (CK,K) \to (C(f),A)$ and $\chi_{g} : (CL,L) \to (C(g),B)$.
Firstly, we define a homeomorphism $\alpha : C([\chi_{f},\chi_{g}]) \approx C(f){\times}C(g)$ by the following.
\begin{align*}&
\alpha\vert_{C(f){\times}B \cup A{\times}C(g)} = \incl : C(f){\times}B \cup A{\times}C(g) \hookrightarrow C(f){\times}C(g),
\\&
\alpha([t,(s{\wedge}x,y)]) = (t{\wedge}x,(ts){\wedge}y),\quad
\alpha([t,(x,s{\wedge}y)]) = ((ts){\wedge}x,t{\wedge}y),
\end{align*}
for any $(s{\wedge}x,y) \in CK{\times}L$, $(x,s{\wedge}y) \in K{\times}CL$ and $t \in [0,1]$
(see Figure \ref{fig-alpha}).
\begin{figure}[!h]
\includegraphics[width=12cm]{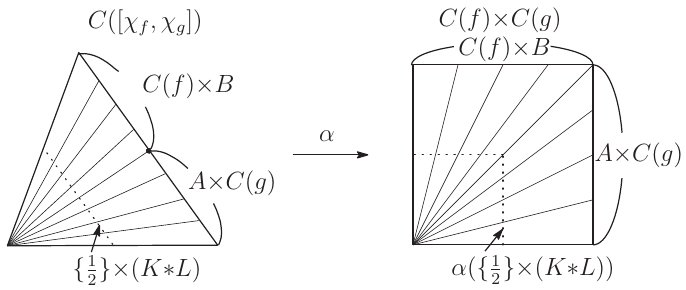}
\caption{}\label{fig-alpha}
\end{figure}\par
Thus the canonical co-pairing $\nu$ is given as follows:
$$
C(f){\times}C(g) \to C(f){\times}C(g)/\alpha(\{C_{{\leq}\frac{1}{2}}(K{\ast}L)\}) \vee C_{{\leq}\frac{1}{2}}(K{\ast}L)/\alpha(\{\frac{1}{2}\}{\times}(K{\ast}L)),
$$
where we can easily see 
\begin{align*}&
C_{{\leq}\frac{1}{2}}(K{\ast}L)/\alpha(\{\frac{1}{2}\}{\times}(K{\ast}L)) 
\\&\qquad
= C_{{\leq}\frac{1}{2}}K{\times}C_{{\leq}\frac{1}{2}}L/((C_{{\leq}\frac{1}{2}}K{\times}\{\frac{1}{2}\}{\times}L) \cup (\{\frac{1}{2}\}{\times}K{\times}C_{{\leq}\frac{1}{2}}L)) \approx \Sigma(K{\ast}L) 
\\[1ex]
&
C(f){\times}C(g)/\alpha(\{C_{{\leq}\frac{1}{2}}(K{\ast}L)\}) 
= C(f){\times}C(g)/C_{\leq\frac{1}{2}}K{\times}C_{\leq\frac{1}{2}}L.
\end{align*}
Secondly, since $C_{{\leq}\frac{1}{2}}X$ is contractible, the inclusions 
$$(C(f),\{\ast\}) \hookrightarrow (C(f),C_{\leq\frac{1}{2}}K) \quad \text{and} \quad (C(g),\{\ast\}) \hookrightarrow (C(g),C_{\leq\frac{1}{2}}L)
$$
are homotopy equivalences, and so is the inclusion
$$
\begin{diagram}
\node{(C(f),\{\ast\}){\times}(C(g),\{\ast\})}
\arrow{e,J}
\node{(C(f),C_{\leq\frac{1}{2}}K){\times}(C(g),C_{\leq\frac{1}{2}}L).}
\end{diagram}
$$
Thus the inclusion
$$
\begin{diagram}
\node{C(f) {\times} \{\ast\} \cup \{\ast\} {\times} C(g)}
\arrow{e,J}
\node{C(f) {\times} C_{\leq\frac{1}{2}}L \cup C_{\leq\frac{1}{2}}K {\times} C(g)}
\end{diagram}
$$
is a homotopy equivalence.
Moreover, since the collapsings
\begin{align*}&
{C(f) {\times} C_{\leq\frac{1}{2}}L \cup C_{\leq\frac{1}{2}}K {\times} C(g)}
\\&\qquad\quad
\longrightarrow
({C(f) {\times} C_{\leq\frac{1}{2}}L \cup C_{\leq\frac{1}{2}}K {\times} C(g))/(C_{\leq\frac{1}{2}}K {\times}C_{\leq\frac{1}{2}}L})\quad \text{and}
\\&
{C(f) {\times} \{\ast\} \cup \{\ast\} {\times} C(g)}
\\&\qquad\quad
\longrightarrow
{(C(f)/C_{\leq\frac{1}{2}}K){\times}\{\ast\} \cup \{\ast\}{\times}(C(g)/C_{\leq\frac{1}{2}}L)}
\end{align*}
are homotopy equivalences, so is the collapsing
\begin{align*}&
({C(f) {\times} C_{\leq\frac{1}{2}}L \cup C_{\leq\frac{1}{2}}K {\times} C(g))/(C_{\leq\frac{1}{2}}K {\times}C_{\leq\frac{1}{2}}L})
\\&\qquad\quad
\longrightarrow
(C(f)/C_{\leq\frac{1}{2}}K){\times}\{\ast\} \cup \{\ast\}{\times}(C(g)/C_{\leq\frac{1}{2}}L)
\\&\qquad\qquad
=
(C_{\geq\frac{1}{2}}(f)/\{\frac{1}{2}\}{\times}K){\times}\{\ast\} \cup \{\ast\}{\times}(C_{\geq\frac{1}{2}}(g)/\{\frac{1}{2}\}{\times}L).
\end{align*}
Finally, since $C_{\leq\frac{1}{2}}K {\times}C_{\leq\frac{1}{2}}L=\alpha(\{C_{{\leq}\frac{1}{2}}(K{\ast}L)\})$, by taking push-out of this collapsing with the inclusion
$$
\divide\dgARROWLENGTH by2
\begin{diagram}
\node{C(f) {\times} C_{\leq\frac{1}{2}}L \cup C_{\leq\frac{1}{2}}K {\times} C(g)/(C_{\leq\frac{1}{2}}K {\times}C_{\leq\frac{1}{2}}L)}
\arrow{s,J}
\\
\node{C(f) {\times} C(g)/\alpha(\{C_{{\leq}\frac{1}{2}}(K{\ast}L)\}),}
\end{diagram}
$$
we obtain a homotopy equivalence
%
\begin{align*}
&
C(f){\times}C(g)/\alpha(\{C_{{\leq}\frac{1}{2}}(K{\ast}L)\}) 
\to C_{{\geq}\frac{1}{2}}(f)/(\{\frac{1}{2}\}{\times}K){\times}C_{{\geq}\frac{1}{2}}(g)/(\{\frac{1}{2}\}{\times}L),
\end{align*}
where the target space is clearly homeomorphic to $C(f){\times}C(g)$.
Therefore, $\nu$ is homotopic to the map $\hat\nu$ given by 
\begin{align*}&
\hat\nu(s{\wedge}x,t{\wedge}y) 
\\&\quad= 
\begin{cases}
(s{\wedge}x,t{\wedge}y) \in C_{\geq\frac{1}{2}}(f)/(\{\frac{1}{2}\}{\times}K) \times C_{\geq\frac{1}{2}}(g)/(\{\frac{1}{2}\}{\times}L),&s, t \geq \frac{1}{2},
\\
(\ast,t{\wedge}y) \in \{\ast\} \times C_{\geq\frac{1}{2}}(g)/(\{\frac{1}{2}\}{\times}L),&s \leq \frac{1}{2}, t \geq \frac{1}{2},
\\
(s{\wedge}x,\ast) \in C_{\geq\frac{1}{2}}(f)/(\{\frac{1}{2}\}{\times}K) \times \{\ast\},&s \geq \frac{1}{2}, t \leq \frac{1}{2},
\\
(s{\wedge}x,t{\wedge}y) \in C_{\leq\frac{1}{2}}K/(\{\frac{1}{2}\}{\times}K) \wedge C_{\leq\frac{1}{2}}L/(\{\frac{1}{2}\}{\times}L),&s, t \leq \frac{1}{2},
\end{cases}
\end{align*}
which coincides exactly the map $\Phi{\circ}(\nu_{f}{\times}\nu_{g})$, and it completes the proof of the proposition.
\end{proof}
\begin{prop}\label{eqnu}
\[
T_1{\circ}((\nu_m {\times} \mathrm{id}_{F'_1}) \vee \mathrm{id}_{\Sigma K^{m,1}_{m+1}})
{\circ}\nu^{m,1}_{m+1}
=
(\nu^{m,1}_{m+1}\cup \mathrm{id}_{\Sigma K_m {\times} F'_1})
 {\circ}(\nu_m {\times} \mathrm{id}_{F'_1}),
\]
where $\nu_m : F_m \to F_m \vee \Sigma K_m$ is the canonical co-pairing
and $T_1: F_{m+1}^{m,1}\cup_{F'_1}(\Sigma K_m {\times} F'_1)\vee\Sigma K^{m,1}_{m+1}
\to (F_{m+1}^{m,1}\vee\Sigma K^{m,1}_{m+1})\cup_{F'_1}(\Sigma K_m {\times} F'_1)$ 
is the canonical homeomorphism.
\end{prop}
\begin{proof}
First, 
%
by Proposition \ref{prop:key-prop}, we obtain the following commutative diagram.
$$
\begin{diagram}
\node{F_{m}{\times}F'_{1}}
	\arrow{e,t}{\nu_{m+1}^{m,1}}
	\arrow{s,l}{\nu_{m}{\times}\id{F'_{1}}}
\node{F_{m}{\times}F_{1} \vee \Sigma(K_{m}{\ast}K'_{1})}
\\
\node{F_{m}{\times}F'_{1} \cup_{F'_{1}}
 \Sigma{K_{m}}{\times}F'_{1}}
	\arrow{e,t}{\id{F_{m}}{\times}\nu_{1}}
\node{(F_{m}{\times}F'_{1} \cup_{F'_{1}}
 \Sigma{K_{m}}{\times}F'_{1} \cup_{F_{m}}
 F_{m}{\times}\Sigma{K'_{1}}) \cup
 \Sigma{K_{m}}{\times}\Sigma{K'_{1}}.}
	\arrow{n,l}{\Phi}
\end{diagram}
$$
Since $\Phi$ goes through $(F_{m}{\times}F'_{1} \cup_{F'_{1}}\Sigma{K_{m}}{\times}F'_{1}) \cup \Sigma{K_{m}}{\times}\Sigma{K'_{1}}/\{\ast\}{\times}\Sigma{K'_{1}}$ as 
$$
\divide\dgARROWLENGTH by2
\begin{diagram}
\node{(F_{m}{\times}F'_{1} \cup_{F'_{1}}
 \Sigma{K_{m}}{\times}F'_{1} \cup_{F_{m}}
 F_{m}{\times}\Sigma{K'_{1}}) \cup
 \Sigma{K_{m}}{\times}\Sigma{K'_{1}}}
\arrow{s,l}{\Phi'}
\\
\node{(F_{m}{\times}F'_{1} \cup_{F'_{1}}
 \Sigma{K_{m}}{\times}F'_{1}) \cup
 \Sigma{K_{m}}{\times}\Sigma{K'_{1}}/\{\ast\}{\times}\Sigma{K'_{1}}}
\arrow{s,l}{p'}
\\
\node{F_{m}{\times}F'_{1} \vee \Sigma(K_{m}{\ast}K'_{1}),}
\end{diagram}
$$
where $\Phi'$ and $p'$ are given by the following formulae.
\begin{align*}&
\Phi'\vert_{F_{m}{\times}F'_{1}} = \id{F_{m}{\times}F'_{1}},
\quad
\Phi'\vert_{\Sigma{K_{m}}{\times}F'_{1}} = \id{\Sigma{K_{m}}{\times}F'_{1}},
\quad
\Phi'\vert_{F_{m}{\times}\Sigma{K'_{1}}} = \proj_{1}
\\&
\Phi'\vert_{\Sigma{K_{m}}{\times}\Sigma{K'_{1}}} = (\text{the canonical collapsing to $\Sigma{K_{m}}{\times}\Sigma{K'_{1}}/\{\ast\}{\times}\Sigma{K'_{1}}$})
\\&
p'\vert\vert_{F_{m}{\times}F'_{1}} = \id{F_{m}{\times}F'_{1}},
\quad
p'\vert_{\Sigma{K_{m}}{\times}F'_{1}} = \proj_{2},
\\&
p'\vert_{\Sigma{K_{m}}{\times}\Sigma{K'_{1}}/\{\ast\}{\times}\Sigma{K'_{1}}} = (\text{the canonical collapsing to $\Sigma(K_{m}{\ast}K'_{1})$})
\end{align*}
which fit in with the following Figure \ref{phi-prime}.
\begin{figure}[!h]
\includegraphics[width=12cm,clip]{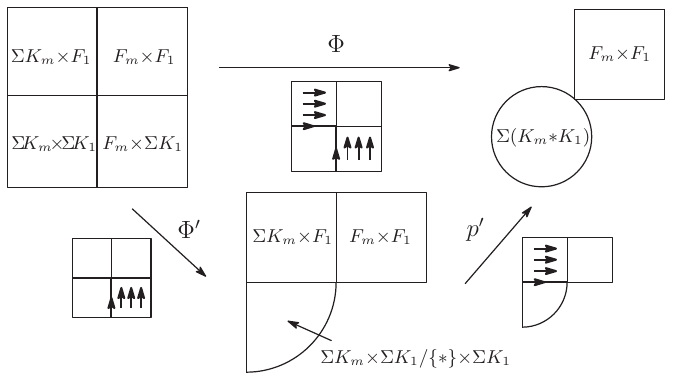}
\caption{}\label{phi-prime}
\end{figure}\par
Since there is a homotopy equivalence $h : \Sigma{K_{m}}{\times}\Sigma{K'_{1}}/\{\ast\}{\times}\Sigma{K'_{1}}$ $\simeq$ $\Sigma{K_{m}}$ $\vee$ $(\Sigma{K_{m}}{\times}\Sigma{K'_{1}}/\Sigma{K_{m}}{\vee}\Sigma{K'_{1}})$ $=$ $\Sigma{K_{m}}$ $\vee$ $\Sigma(K_{m}{\ast}K'_{1})$ such that $h\vert_{\Sigma{K_{m}}{\times}\{\ast\}}=\id{\Sigma{K_{m}}}$, we obtain that $p'$ goes through $(F_{m}{\times}F'_{1} \cup_{F'_{1}} \Sigma{K_{m}}{\times}F'_{1}) \vee \Sigma(K_{m}{\ast}K'_{1})$ as
$$
\divide\dgARROWLENGTH by2
\begin{diagram}
\node{(F_{m}{\times}F'_{1} \cup_{F'_{1}}
 \Sigma{K_{m}}{\times}F'_{1}) \cup
 \Sigma{K_{m}}{\times}\Sigma{K'_{1}}/\{\ast\}{\times}\Sigma{K'_{1}}}
\arrow{s,l}{p_{0}}
\\
\node{(F_{m}{\times}F'_{1} \cup_{F'_{1}}
 \Sigma{K_{m}}{\times}F'_{1}) \vee \Sigma(K_{m}{\ast}K'_{1})}
\arrow{s,l}{p_{1}}
\\
\node{F_{m}{\times}F'_{1} \vee \Sigma(K_{m}{\ast}K'_{1}),}
\end{diagram}
$$
where $p_{0}$ and $p_{1}$ are given by the following formulae.
\begin{align*}&
p_{0}\vert_{F_{m}{\times}F'_{1}} = \id{F_{m}{\times}F'_{1}},\quad
p_{0}\vert_{\Sigma{K_{m}}{\times}F'_{1}} = \id{\Sigma{K_{m}}{\times}F'_{1}},\quad
p_{0}\vert_{\Sigma{K_{m}}{\times}\Sigma{K'_{1}}/\{\ast\}{\times}\Sigma{K'_{1}}} = h
\\&
p_{1}\vert_{F_{m}{\times}F'_{1}} = \id{F_{m}{\times}F'_{1}},\quad
p_{1}\vert_{\Sigma{K_{m}}{\times}F'_{1}} = \proj_{2},\quad
p_{1}\vert_{\Sigma(K_{m}{\ast}K'_{1})} = \id{\Sigma(K_{m}{\ast}K'_{1})}
\end{align*}
which fit in with the following Figure \ref{p_0-p_1}.%
\begin{figure}[!h]
\includegraphics[width=12cm,clip]{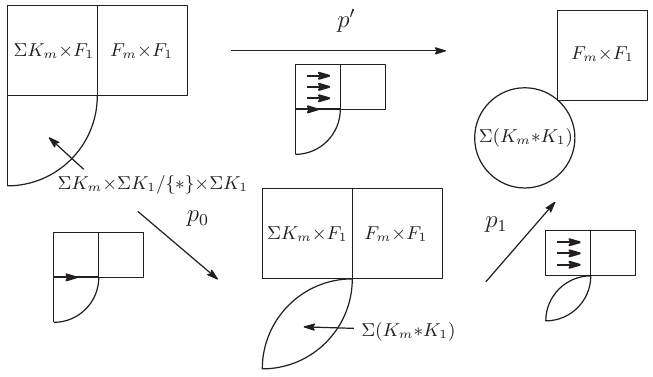}
\caption{}\label{p_0-p_1}
\end{figure}\par
Hence $\Phi$ is decomposed as $\Phi=p'{\circ}\Phi'=p_{1}{\circ}p_{0}{\circ}\Phi'$ and $p_{0}{\circ}\Phi'$ can be written as
\begin{align*}&
p_{0}{\circ}\Phi'\vert_{F_{m}{\times}F'_{1}} = \id{F_{m}{\times}F'_{1}},
\quad
p_{0}{\circ}\Phi'\vert_{\Sigma{K_{m}}{\times}F'_{1}} = \id{\Sigma{K_{m}}{\times}F'_{1}},
\quad
p_{0}{\circ}\Phi'\vert_{F_{m}{\times}\Sigma{K'_{1}}} = \proj_{1}
\\&
p_{0}{\circ}\Phi'\vert_{\Sigma{K_{m}}{\times}\Sigma{K'_{1}}} : \Sigma{K_{m}}{\times}\Sigma{K'_{1}} \to \Sigma{K_{m}} \vee \Sigma(K_{m}{\ast}K'_{1})
\end{align*}
and hence $p_{0}{\circ}\Phi'{\circ}(\id{F_{m}}{\times}\nu_{1})$ is given by
\begin{align*}&
p_{0}{\circ}\Phi'{\circ}(\id{F_{m}}{\times}\nu_{1})\vert_{F_{m}{\times}F'_{1}} = \id{F_{m}{\times}F'_{1}},
\\&
p_{0}{\circ}\Phi'{\circ}(\id{F_{m}}{\times}\nu_{1})\vert_{\Sigma{K_{m}}{\times}F'_{1}}=\nu^{\prime} : \Sigma{K_{m}}{\times}F'_{1} \to \Sigma{K_{m}}{\times}F'_{1} \vee \Sigma(K_{m}{\ast}K'_{1}),
\end{align*}
where $\nu^{\prime}$ is giving a co-pairing.
Thus we obtain a commutative diagram
\begin{equation}\label{firsteq}
\xymatrix{
\hspace{-15mm}F_{m+1}^{m,1}=F_{m}{\times}F'_{1} \ar[r]^{\nu_m{\times} {\rm id}_{F'_1}\hspace{30pt}} \ar[d]^{\nu^{m,1}_{m+1}}
 & F_{m}{\times}F'_{1}\cup_{F'_1} (\Sigma K_m {\times} F'_1)
     \ar[d]^{{\rm id}_{F_{m}{\times}F'_{1}}\cup \nu^{\prime}} \\
F_{m}{\times}F'_{1} \vee \Sigma K_m {\ast} K'_1
 & F_{m}{\times}F'_{1}\cup_{F'_1} (\Sigma K_m {\times} F'_1) \vee \Sigma K_m {\ast} K'_1.
\ar[l]_{p_1\hspace{30pt}}
}
\end{equation}
%
%
Therefore we have
\[
\begin{split}
T_1{\circ}((\nu_m {\times}& {\rm id}_{F'_1}) \vee {\rm id}_{\Sigma K^{m,1}_{m+1}})
{\circ}\nu^{m,1}_{m+1}\\
&=
T_1{\circ}((\nu_m {\times} {\rm id}_{F'_1}) \vee {\rm id}_{\Sigma K^{m,1}_{m+1}})
 {\circ}p_1{\circ}({\rm id}_{F_{m+1}^{m,1}}\cup \nu^{\prime})
 {\circ}(\nu_m{\times} {\rm id}_{F'_1}).
\end{split}
\]
Let us denote
$
p_2 : F_{m+1}^{m,1} \cup_{F'_1} (\Sigma K_m {\times} F'_1) 
        \cup_{F'_1} (\Sigma K_m {\times} F'_1)\vee \Sigma K^{m,1}_{m+1}
 \to F_{m+1}^{m,1} \cup_{F'_1} (\Sigma K_m {\times} F'_1)\vee \Sigma K^{m,1}_{m+1}
$
by the map pinching the second $\Sigma K_m {\times} F'_1$ to $F'_{1}$,
$
p_3 : F_{m+1}^{m,1}\cup_{F'_1}((\Sigma K_m {\times} F'_1)\vee \Sigma K^{m,1}_{m+1})
      \cup_{F'_1}(\Sigma K_m {\times} F'_1)
 \to (F_{m+1}^{m,1}\vee\Sigma K^{m,1}_{m+1})\cup_{F'_1}\Sigma K^{m,1}_{m+1}
$
by the map pinching the first $\Sigma K_m {\times} F'_1$ to one point,
$\nu_0 : \Sigma K_m \to \Sigma K_m \vee \Sigma K_m$ 
by the canonical co-multiplication and 
$T_0 :\Sigma K_m \vee \Sigma K_m \to \Sigma K_m \vee \Sigma K_m$ by the commutative map.
It is easy to check the following.
\[
\begin{split}&
T_1{\circ}((\nu_m {\times} {\rm id}_{F'_1}) \vee {\rm id}_{\Sigma K^{m,1}_{m+1}})
 {\circ}\nu^{m,1}_{m+1}
\\&\qquad=
T_1{\circ}p_2{\circ}((\nu_m{\times}{\rm id}_{F'_1})
     \cup{\rm id}_{\Sigma K_m {\times} F'_1}\vee {\rm id}_{\Sigma K_m {\ast} K'_1})
{\circ}({\rm id}_{F_{m+1}^{m,1}}\cup \nu^{\prime})
   {\circ}(\nu_m{\times} {\rm id}_{F'_1})
\\&\qquad=
T_1{\circ}p_2 
 {\circ} ({\rm id}_{F_{m+1}^{m,1}}\cup {\rm id}_{\Sigma K_m {\times} F'_1}
         \cup \nu^{\prime})
{\circ}((\nu_m{\times} {\rm id}_{F'_1})\cup {\rm id}_{\Sigma K_m {\times} F'_1})
 {\circ}(\nu_m{\times} {\rm id}_{F'_1})
\\&\qquad=
p_3{\circ}({\rm id}_{F_{m+1}^{m,1}}\cup\nu^{\prime}\cup {\rm id}_{\Sigma K_m {\times} F'_1})
  {\circ}({\rm id}_{F_{m+1}^{m,1}}\cup(T_0{\times} {\rm id}_{F'_1}))
\\&\qquad\hspace{160pt}
{\circ}((\nu_m{\times} {\rm id}_{F'_1})\cup {\rm id}_{\Sigma K_m {\times} F'_1})
 {\circ}(\nu_m{\times} {\rm id}_{F'_1}).
\end{split}
\]
Using the equations 
$({\rm id}_{F_m}{\times} \nu_0){\circ}\nu_m
 = (\nu_m {\times} {\rm id}_{F_m}){\circ}\nu_m$
and $T_0{\circ}\nu_0 = \nu_0$ from the assumption that $K_m$ is a co-H-space,
we have 
\[
\begin{split}
&
T_1{\circ}((\nu_m {\times} {\rm id}_{F'_1}) \vee {\rm id}_{\Sigma K^{m,1}_{m+1}})
 {\circ}\nu^{m,1}_{m+1}
\\&\qquad=
p_3{\circ}({\rm id}_{F_{m+1}^{m,1}}\cup\nu^{\prime}\cup {\rm id}_{\Sigma K_m {\times} F'_1})
 {\circ}({\rm id}_{F_{m+1}^{m,1}}\cup(T_0{\times} {\rm id}_{F'_1}))
\\&\qquad\hspace{160pt}
{\circ}({\rm id}_{F_{m+1}^{m,1}}\cup(\nu_0 {\times} {\rm id}_{F'_1}))
 {\circ}(\nu_m{\times} {\rm id}_{F'_1})
\\&\qquad=
p_3{\circ}({\rm id}_{F_{m+1}^{m,1}}\cup\nu^{\prime}\cup {\rm id}_{\Sigma K_m {\times} F'_1})
{\circ}({\rm id}_{F_{m+1}^{m,1}}\cup(\nu_0 {\times} {\rm id}_{F'_1}))
 {\circ}(\nu_m{\times} {\rm id}_{F'_1})
\\&\qquad=
p_3{\circ}({\rm id}_{F_{m+1}^{m,1}}\cup\nu^{\prime}\cup {\rm id}_{\Sigma K_m {\times} F'_1})
{\circ}((\nu_m{\times} {\rm id}_{F'_1})\cup {\rm id}_{\Sigma K_m {\times} F'_1})
 {\circ}(\nu_m{\times} {\rm id}_{F'_1}).
\end{split}
\]
Using the diagram (\ref{firsteq}), we proceed further as follows:
\[
T_1{\circ}((\nu_m {\times} {\rm id}_{F'_1}) \vee {\rm id}_{\Sigma K^{m,1}_{m+1}})
 {\circ}\nu^{m,1}_{m+1}
=
(\nu_{m+1}^{m,1}\cup {\rm id}_{\Sigma K_m {\times} F'_1})
 {\circ}(\nu_m{\times} {\rm id}_{F'_1}).
\]
This completes the proof of Proposition \ref{eqnu}.
\end{proof}
\begin{proof}[{\it Proof of Lemma \ref{2cofib}}]
The commutativity of the left-hand side square follows from Proposition 2.9 of \cite{Stanley02} and the middle square is clearly commutative.
To see the formula $(\sigma_m{\times} \sigma'\vee \Sigma g_m{\ast} g'){\circ}\nu^{m,1}_{m+1}  = \hat{\nu}_{m+1}{\circ}(\sigma_m{\times} \sigma')$, let us recall the construction of the structure map $\sigma_m : F_m \to P^{m}\Omega F_m$.
Then we obtain $\sigma_m = \nabla_{P^{m}\Omega F_m} {\circ}(\sigma^{\prime}_m \vee \iota^{\Omega F_m}_{m-1,m}{\circ}\delta_m){\circ}\nu_m$, where $\sigma^{\prime}_m$ is induced from the diagram
$$\xymatrix{
K_m \ar[r]^{h_{m}} \ar[d]^{g_{m}^{\prime}}
 & F_{m-1} \ar@{^{(}->}[rr]^{i_{m-1,m}^F} \ar[d]^{P^{m-1}\Omega i^{F}_{m-1,m}{\circ}\sigma_{m-1}}
 &
 & F_{m}  \ar@{..>}[d]^{\sigma^{\prime}_{m}}
\\
E^{m}\Omega F_{m} \ar[r]_{P_{m}^{\Omega F_m}}
 & P^{m-1}\Omega F_{m} \ar@{^{(}->}[rr]_{\iota_{m-1,m}^{\Omega F_m}}
 &
 & P^{m}\Omega F_{m}.
}$$
The map $\delta_m : \Sigma K_m \to P^{m-1}{\Omega F_{m}}$ is the map pulled back the difference map $\delta^{\prime}_m : \Sigma K_m \to F_m$ which is the difference between the identity map of $F_m$ and $e^{F_m}_m{\circ}\sigma_m^{\prime}$.
Thus we have the equation
\begin{align*}&
(\sigma_m{\times} \sigma' \vee \Sigma g_m{\ast} g'){\circ}\nu^{m,1}_{m+1}
\\&\quad
 = \{(\nabla_{P^{m}\Omega F_m}{\circ}
     (\sigma^{\prime}_m \vee (\iota^{\Omega F_m}_{m-1,m}{\circ}\delta_m)){\circ}\nu_m){\times}
     \sigma'\vee \Sigma g_m{\ast} g' \}{\circ}\nu^{m,1}_{m+1}
\\&\quad
 = \{(\nabla_{P^{m}\Omega F_m} {\times} {\rm id}_{\Sigma\Omega F'_1}){\circ}((
     \sigma_m^{\prime} \vee (\iota^{\Omega F_m}_{m-1,m}{\circ}\delta_m)){\times} \sigma') {\circ}(\nu_m{\times} {\rm id}_{F'_1})\vee \Sigma g_m{\ast} g'\}
\\&\quad \hspace{30pt}
 {\circ}\nu^{m,1}_{m+1}
\\&\quad
 = (\nabla_{P^{m}\Omega F_m} {\times} {\rm id}_{\Sigma\Omega F'_1} \vee {\rm id}_{\Sigma\hat{E}_{m+1}})
    {\circ}\{((\sigma_m^{\prime} \vee (\iota^{\Omega F_m}_{m-1,m}{\circ}\delta_m)){\times} \sigma')
     \vee \Sigma g_m{\ast} g'\}
\\&\quad \hspace{30pt}
 {\circ}((\nu_m {\times} {\rm id}_{F'_1}) \vee {\rm id}_{\Sigma K^{m,1}_{m+1}}){\circ}\nu^{m,1}_{m+1}
\\&\quad
 = (\nabla_{P^{m}\Omega F_m} {\times} {\rm id}_{\Sigma\Omega F'_1} \vee {\rm id}_{\Sigma\hat{E}_{m+1}})
    {\circ}\{(\sigma_m^{\prime}{\times} \sigma')
     \cup ((\iota^{\Omega F_m}_{m-1,m}{\circ}\delta_m){\times} \sigma')
     \vee \Sigma g_m{\ast} g'\}
\\&\quad \hspace{30pt}
   {\circ}((\nu_m {\times} {\rm id}_{F'_1}) \vee {\rm id}_{\Sigma K^{m,1}_{m+1}}){\circ}\nu^{m,1}_{m+1}
\\&\quad
 = (\nabla_{P^{m}\Omega F_m} {\times} {\rm id}_{\Sigma\Omega F'_1} \vee {\rm id}_{\Sigma\hat{E}_{m+1}})
\\&\quad \hspace{30pt}
    {\circ}T_2{\circ}\{(\sigma_m^{\prime}{\times} \sigma' \vee \Sigma g_m{\ast} g')
     \cup ((\iota^{\Omega F_m}_{m-1,m}{\circ}\delta_m){\times} \sigma')\}
\\&\quad \hspace{30pt}
   {\circ}T_1{\circ}((\nu_m {\times} {\rm id}_{F'_1}) \vee {\rm id}_{\Sigma K^{m,1}_{m+1}}){\circ}\nu^{m,1}_{m+1},
\end{align*}
where $T_2 : (\hat{F}_{m+1}\vee \Sigma \hat{E}_{m+1})\cup_{\Sigma\Omega F'_1} \hat{F}_{m+1}
\to (\hat{F}_{m+1}\cup_{\Sigma\Omega F'_1} \hat{F}_{m+1})\vee \Sigma \hat{E}_{m+1}$
is an appropriate homeomorphism.
By Proposition \ref{eqnu}, we can proceed as follows:
\begin{align*}&
(\sigma_m{\times}  \sigma' \vee \Sigma g_m{\ast} g'){\circ}\nu^{m,1}_{m+1}
  = (\nabla_{P^{m}\Omega F_m} {\times} {\rm id}_{\Sigma\Omega F'_1} \vee {\rm id}_{\Sigma\hat{E}_{m+1}})
\\
& \hspace{160pt}
    {\circ}T_2{\circ}\{(\sigma_m^{\prime}{\times} \sigma' \vee \Sigma g_m{\ast} g')
     \cup ((\iota^{\Omega F_m}_{m-1,m}{\circ}\delta_m){\times} \sigma')\}
\\ 
& \hspace{160pt}
    {\circ}(\nu^{m,1}_{m+1}\cup {\rm id}_{\Sigma K_m {\times} F'_1}){\circ}(\nu_m {\times} {\rm id}_{F'_1})
\\&\qquad
  = (\nabla_{P^{m}\Omega F_m} {\times} {\rm id}_{\Sigma\Omega F'_1} \vee {\rm id}_{\Sigma\hat{E}_{m+1}})
    {\circ}T_2
\\&\qquad \hspace{30pt}
    {\circ}\{((\sigma_m^{\prime}{\times} \sigma' \vee \Sigma g_m{\ast} g'){\circ}\nu^{m,1}_{m+1})
     \cup ((\iota^{\Omega F_m}_{m-1,m}{\circ}\delta_m){\times} \sigma')\}
   {\circ}(\nu_m {\times} {\rm id}_{F'_1}).
\end{align*}

By the definitions of $\sigma_m^{\prime}$ and $\sigma'$, we have
\begin{align*}
(\sigma_m{\times} & \sigma' \vee \Sigma g_m{\ast} g'){\circ}\nu^{m,1}_{m+1}
\\
& = (\nabla_{P^{m}\Omega F_m} {\times} {\rm id}_{\Sigma\Omega F'_1} \vee {\rm id}_{\Sigma\hat{E}_{m+1}})
   {\circ}T_2
\\ 
& \hspace{30pt}
    {\circ}\{(\hat{\nu}_{m+1}{\circ}(\sigma_m^{\prime}{\times} \sigma'))
     \cup ((\iota^{\Omega F_m}_{m-1,m}{\circ}\delta_m){\times} \sigma')\}
   {\circ}(\nu_m {\times} {\rm id}_{F'_1})
\\
& = (\nabla_{P^{m}\Omega F_m} {\times} {\rm id}_{\Sigma\Omega F'_1} \vee \nabla_{\Sigma\hat{E}_{m+1}})
   {\circ}T_3
\\ 
& \hspace{30pt}
    {\circ}\{\hat{\nu}_{m+1}{\circ}(\sigma_m^{\prime}{\times} \sigma')
     \cup i_1{\circ}((\iota^{\Omega F_m}_{m-1,m}{\circ}\delta_m){\times} \sigma')\}
%
   {\circ}(\nu_m {\times} {\rm id}_{F'_1}).
\end{align*}
Here $i_1 : \hat{F}_{m+1} \to \hat{F}_{m+1} \vee \Sigma \hat{E}_{m+1}$ is the inclusion map and 
$T_3 : (\hat{F}_{m+1}\vee \Sigma \hat{E}_{m+1})\cup_{\Sigma\Omega F'_1} (\hat{F}_{m+1}\vee \Sigma \hat{E}_{m+1}) \to (\hat{F}_{m+1}\cup_{\Sigma\Omega F'_1} \hat{F}_{m+1})\vee \Sigma \hat{E}_{m+1} \vee \Sigma \hat{E}_{m+1}$ is the canonical homeomorphism.
\begin{align*}
(\sigma_m{\times} & \sigma' \vee \Sigma g_m{\ast} g'){\circ}\nu^{m,1}_{m+1}
\\
& = (\nabla_{P^{m}\Omega F_m} {\times} {\rm id}_{\Sigma\Omega F'_1} \vee \nabla_{\Sigma\hat{E}_{m+1}})
  {\circ}T_3{\circ}(\hat{\nu}_{m+1}\cup \hat{\nu}_{m+1} ) 
\\ 
& \hspace{30pt}
    {\circ}\{(\sigma_m^{\prime}{\times} \sigma')
     \cup ((\iota^{\Omega F_m}_{m-1,m}{\circ}\delta_m){\times} \sigma')\}
   {\circ}(\nu_m {\times} {\rm id}_{F'_1})
\\
& = \hat{\nu}_{m+1} 
   {\circ}(\nabla_{P^{m}\Omega F_m} {\times} {\rm id}_{\Sigma\Omega F'_1})
   {\circ}\{(\sigma_m^{\prime}{\times} \sigma')
    \cup ((\iota^{\Omega F_m}_{m-1,m}{\circ}\delta_m){\times} \sigma')\}
   {\circ}(\nu_m {\times} {\rm id}_{F'_1})
\\
& = \hat{\nu}_{m+1} 
   {\circ}(\nabla_{P^{m}\Omega F_m} {\times} {\rm id}_{\Sigma\Omega F'_1})
   {\circ}\{(\sigma_m^{\prime}\vee (\iota^{\Omega F_m}_{m-1,m}{\circ}\delta_m)) {\times} \sigma'\}
   {\circ}(\nu_m {\times} {\rm id}_{F'_1})
\\
& = \hat{\nu}_{m+1} 
   {\circ}\{\nabla_{P^{m}\Omega F_m}
   {\circ}(\sigma_m^{\prime}\vee (\iota^{\Omega F_m}_{m-1,m}{\circ}\delta_m))
   {\circ}\nu_m  {\times} \sigma'\}
 = \hat{\nu}_{m+1}{\circ}(\sigma_m{\times}\sigma').
\end{align*}
This completes the proof.
\end{proof}
\section{Proof of Theorem \ref{main2thm}}%
\label{proof-mthm}
In the fibre sequence 
$G\hookrightarrow E\to \Sigma A$,
by the James-Whitehead decomposition
(see Theorem VII.(1.15) of Whitehead \cite{GWhitehead78}),
 the total space $E$ has the homotopy type of the space
$G\cup_{\psi}G{\times} CA$.
Here $\psi$ is the following composition.
\[
\psi : G{\times} A\xrightarrow{\mathrm{id}_{G} {\times} \alpha} G{\times} G\xrightarrow{\mu} G.
\]
Since $G \simeq F_m$ and $\alpha$ is compressible into $F^{\prime}_1$,
we can see that
\[
\psi : G {\times} A \simeq F_m{\times} A \xrightarrow{\mathrm{id}_{F_m}{\times} \alpha}
 F_m {\times} F^{\prime}_1 \subset F_m {\times} F_1 \subset F_m {\times} F_m \simeq G{\times} G
  \xrightarrow{\mu} G \simeq F_m
\]
and $E$ is the homotopy push out of the following sequence.
\[
\xymatrix{
F_m 
& & F_m {\times} A \ar[ll]_{pr_1} \ar[rr]^{\mathrm{id}_{F_m}{\times} \alpha}
& & F_m {\times} F'_1 \ar[rr]^{\hspace{15pt}\mu'_{m,1}}
& & F_m.
}\]
We construct spaces and maps
such that the homotopy push out of these maps dominates $E$.
\par
The condition of $H_1(\alpha )=0 \in [A, \Omega F'_1{\ast} \Omega F'_1]$
implies that 
\begin{equation}
\label{Ead-a}
\Sigma \ad{\alpha} 
= \sigma'{\circ}\alpha : A \to F^{\prime}_{1}\to \Sigma \Omega F'_{1}.
\end{equation}%
Let $e' : \Omega\Sigma F'_1\to F'_1$ and $e^A_1 : \Sigma\Omega A \to A$ be evaluation maps. 
Since $A$ is a suspended space, 
we put $\sigma_A= \Sigma\ad{\mathrm{id}_A} : A \to \Sigma\Omega A$. 
\begin{lem}\label{diagPxP}
The following diagram is commutative.
\[\xymatrix{
F_m \ar[d]^{\iota^{\Omega F_m}_{m,m+1}{\circ}\sigma_m}
 & F_m{\times} A \ar[d]^{\sigma_m{\times} \sigma_A} \ar[l]_{pr_1} \ar[r]^{\mathrm{id}_{F_m}{\times} \alpha }
 & F_m{\times} F'_1 \ar[d]^{\sigma_m{\times} \sigma'}  \ar[r]^{\mu'_{m,1}}
 & F_m \ar[d]^{\iota^{\Omega F_m}_{m,m+1}{\circ}\sigma_m}\\
P^{m+1}\Omega F_m \ar[d]^{e^{F_m}_{m+1}}
 & (P^m\Omega F_m)^{(\ell)}{\times} (\Sigma \Omega A)^{(\ell)}
    \ar[l]_{\hspace{-20pt}\phi}
    \ar[d]^{(e^{F_m}_m)^{(\ell)}{\times} (e^A_1)^{(\ell)}}
    \ar[r]^{\hspace{40pt}\chi}
 & \hat{F}_{m+1} \ar[d]^{(e^{F_m}_m)^{(\ell)}{\times} (e')^{(\ell)}}
 & P^{m+1}\Omega F_m \ar[d]^{e^{F_m}_{m+1}}\\
F_m
 & F_m{\times} A \ar[l]_{pr_1} \ar[r]^{\mathrm{id}_{F_m}{\times} \alpha }
 & F_m{\times} F'_1 \ar[r]^{\mu'_{m,1}}
 & F_m,
 }\]%
where $\phi = (\iota^{\Omega F_m}_{m,m+1})^{(\ell)}{\circ}pr_1$
and $\chi = \mathrm{id}_{(P^m\Omega F_m)^{(\ell)}} {\times} (\Sigma \Omega \alpha)^{(\ell)}$.
\end{lem}
\begin{proof}
It is clear to see that the left upper square is commutative.
The equation $e^{F_m}_m = e^{F_m}_{m+1}{\circ}\iota^{\Omega F_m}_{m,m+1}$ implies that the left lower square is commutative.
The commutativity of the middle lower square follows from the equation $\alpha{\circ}e^{A}_{1} = e'{\circ}\Sigma \Omega \alpha$.
Equation (\ref{Ead-a}) implies that the middle upper square is commutative.
Since $\sigma_m$ satisfy the condition (2) of Proposition \ref{cd-to-gs} and  $e'{\circ}\sigma'=\mathrm{id}_{F'_1}$, the right rectangular is commutative, too. This completes the proof of the lemma.
\end{proof}
Let $\lambda=\mu'_{m,1}{\circ}\{(e^{F_m}_m)^{(\ell)}{\times} (e')^{(\ell)}\} : \hat{F}_{m+1} \to F_{m}{\times} F'_1 \to F_m$.
Then $\lambda$ is a well-defined filtered map w.r.t. the filtration (\ref{cd-pxp}) of $\hat{F}_{m+1}$ and the trivial filtration ($(F_{m})_{i}=F_{m}$ for all $i$) of $F_{m}$, where we have $(e^{F_m}_m)^{(\ell)}{\times} (e')^{(\ell)}(\hat{F}_{k}) = {(e^{F_{m-1}}_k)^{(\ell)}{\times} \{\ast\}} \cup (e^{F_{m-1}}_{k-1})^{(\ell)}{\times} (e')^{(\ell)}(\hat{F}_{k}) \subset F_{m-1}{\times} \{\ast\} \cup F_{m-1}{\times} F'_1$ for $0 \leq k \leq m{-}1$ and $(e^{F_m}_m)^{(\ell)}{\times} (e')^{(\ell)}(\hat{F}_{m}) = {(e^{F_{m}}_m)^{(\ell)}{\times} \{\ast\}} \cup (e^{F_{m-1}}_{m-1})^{(\ell)}{\times} (e')^{(\ell)}(\hat{F}_{m}) \subset F_{m}{\times} \{\ast\} \cup F_{m-1}{\times} F'_1$.
\begin{defi}
By applying Proposition \ref{cd-to-gs} to $\lambda$, we obtain a series of maps
\[
\hat{\lambda}_k : \hat{F}_{k}\to P^k\Omega F_m, \qquad 0\le k\le m{+}1.
\]
\end{defi}
\begin{lem}\label{hatmu}
There is a map $\hat{\lambda} : \hat{F}_{m+1} \to P^{m+1}\Omega F_m$
which fits into the right-hand commutative square of the diagram given in Lemma \ref{diagPxP} to make the resulting upper and lower squares sharing $\hat\lambda$ commutative.
\end{lem}
\begin{proof}
%
\par
First, for $0\leq k\leq m{-}1$, we assert 
\begin{equation}\label{su=us}
\iota^{\Omega F_m}_{k,k+1}{\circ}P^k\Omega i^F_{k,m}{\circ}\sigma_k{\circ}\mu^{m,1}_k
 = \iota^{\Omega F_m}_{k,k+1}{\circ}\hat{\lambda}_k{\circ}j_k
  {\circ}\hat{\sigma}_k : F^{m,1}_k \to P^{k+1}\Omega F_m,
\end{equation}
where
$\mu^{m,1}_k=\mu_{k,0}\cup\mu'_{k-1,1} : F^{m,1}_k= F_k{\times} \{\ast\} \cup F_{k-1}{\times} F'_1 \to F_k$,
\[
\hat{\sigma}_k = \sigma_k {\times} \{\ast\}\cup \sigma_{k-1}{\times} \sigma' : F^{m,1'}_k \to
 (P^{k}\Omega F_{k})^{(\ell)}{\times} \{\ast\} \cup (P^{k-1}\Omega F_{k-1})^{(\ell)}{\times} (\Sigma\Omega F'_1)^{(\ell)}
\]
and $j_k=(P^{k}\Omega i_{k,m-1}^{F})^{(\ell)}{\times} \{\ast\}
 \cup (P^{k-1}\Omega i_{k-1,m-1}^{F})^{(\ell)}{\times} \mathrm{id}_{(\Sigma\Omega F'_1)^{(\ell)}}$.
Note that this condition is natural to cone-decompositions.

We show (\ref{su=us}) by induction on $k$.
The case $k=0$ is clear,
since both maps are constant maps.
Assume
the $k$-th of (\ref{su=us}).
Since $\sigma_i$ satisfy the condition (1) of Proposition \ref{cd-to-gs}, 
the following diagram is commutative 
\[\xymatrix{
F_{i} \ar[r]^{\sigma_i} \ar[d]^{i^{F}_{i,i+1}}
 & P^{i}\Omega F_{i} \ar[r]^{P^{i}\Omega i^{F}_{i,i+1}}
 & P^{i}\Omega F_{i+1} \ar@{^{(}->}[rr]^{P^{i}\Omega i^{F}_{i+1,m-1}} \ar@{^{(}->}[d]^{\iota^{\Omega F_{i+1}}_{i,i+1}}
 &
 & P^{i}\Omega F_{m-1} \ar@{^{(}->}[d]^{\iota^{\Omega F_{m-1}}_{i,i+1}}
\\
F_{i+1} \ar[rr]^{\sigma_{i+1}}
 & 
 & P^{i+1}\Omega F_{i+1} \ar@{^{(}->}[rr]^{P^{i+1}\Omega i^{F}_{i+1,m-1}}
 &
 & P^{i+1}\Omega F_{m-1} 
}\]
for $1\le i\le m-2$.
Hence, we have
\begin{align*}
j_{k+1}&{\circ}\hat{\sigma}_{k+1}{\circ}i^{m,1}_{k}\\
& =((P^{k+1}\Omega i_{k+1,m-1}^{F})^{(\ell)}{\circ}\sigma_{k+1}{\circ}i^F_{k,k+1}){\times}\{\ast\}\cup((P^{k}\Omega i_{k,m-1}^{F})^{(\ell)}{\circ}\sigma_k{\circ}i^F_{k-1,k}){\times}\sigma'\\
& =((\iota^{\Omega F_{m-1}}_{k,k+1}{\circ}P^{k}\Omega i^{F}_{k,m-1})^{(\ell)}{\circ}\sigma_k){\times}\{\ast\}\cup ((\iota^{\Omega F_{m-1}}_{k-1,k}{\circ}P^{k}\Omega i^{F}_{k-1,m-1})^{(\ell)}{\circ}\sigma_{k-1}){\times}\sigma'\\
& = \hat{i}_k{\circ}j_k{\circ}\sigma_k.
\end{align*}
By the condition (1) of Proposition \ref{cd-to-gs} for $\hat{\lambda}_{i}$,
we obtain 
$
\hat{\lambda}_{k+1}{\circ}\hat{i}_k = \iota^{\Omega F_m}_{k,k+1}\circ\hat{\lambda}_{k}.
$
Thus we have the equation
\begin{align*}&
{i^{m,1}_{k}}^{\ast}(
 \hat{\lambda}_{k+1}{\circ}j_{k+1}{\circ}\hat{\sigma}_{k+1})
  = 
      \hat{\lambda}_{k+1}{\circ}j_{k+1}{\circ}\hat{\sigma}_{k+1}{\circ}i^{m,1}_{k}
\\&\qquad
  = 
    \hat{\lambda}_{k+1}{\circ}\hat{i}_k{\circ}j_k{\circ}\hat{\sigma}_{k}
  = \iota^{\Omega F_m}_{k,k+1}{\circ}\hat{\lambda}_{k}{\circ}j_k{\circ}\hat{\sigma}_{k}.
\intertext{By the condition (1) of Proposition \ref{cd-to-gs} for $\sigma_{k}$ and the induction hypothesis, we proceed further as follows:}
&\qquad
  = \iota^{\Omega F_m}_{k,k+1}{\circ}P^k\Omega i^F_{k,m}{\circ}\sigma_k{\circ}\mu^{m,1}_k
  = P^{k+1}\Omega i^F_{k+1,m}{\circ}\sigma_{k+1}{\circ}i^F_{k,k+1}{\circ}\mu^{m,1}_k
\\&\qquad
  = {i^{m,1}_{k}}^{\ast}(P^{k+1}\Omega i^F_{k+1,m}{\circ}\sigma_{k+1}{\circ}\mu^{m,1}_{k+1}).
\end{align*}
Thus by a standard argument of homotopy theory, we obtain the difference map
$\delta_{k+1} : \Sigma K^{m,1}_{k+1}\to P^{k+1}\Omega F_m$
such that 
\begin{equation}\label{delta_k+1}
P^{k+1}\Omega i^F_{k+1,m}{\circ}\sigma_{k+1}{\circ}\mu^{m,1}_{k+1} = \nabla_{P^{k+1}\Omega F_m}
{\circ}(\hat{\lambda}_{k+1}{\circ}j_{k+1}{\circ}\hat{\sigma}_{k+1}
\vee \delta_{k+1}){\circ}\nu^{m,1}_{k+1}.
\end{equation}
By the condition (2) of Proposition \ref{cd-to-gs} of $\hat{\lambda}_{k+1}$,
we have the equation
\[
 e_{k+1}^{F_m}{\circ}\hat{\lambda}_{k+1}
  = (\mu_{m-1,0}\cup\mu'_{m-1,1}){\circ}\{{(e^{F_{m-1}}_{k+1})^{(\ell)}{\times} \{\ast\}}
   \cup (e^{F_{m-1}}_{k})^{(\ell)}{\times} (e')^{(\ell)}\}.
\]
By the commutative diagram
\[\xymatrix{
F_i \ar[r]^{\hspace{-10pt}\sigma_i} \ar[dr]_{\mathrm{id}_{F_i}}
 & (P^{i}\Omega F_{i})^{(\ell)} \ar@{^{(}->}[rr]^{(P^{i}\Omega i^{F}_{i,m-1})^{(\ell)}} \ar[d]^{(e^{F_i}_i)^{(\ell)}}
 &
 & (P^{i}\Omega F_{m-1})^{(\ell)} \ar[rr]^{\hspace{20pt}(e^{F_{m-1}}_i)^{(\ell)}}
 &
 & F_{m-1}\\
 & F_i \ar@{^{(}->}[rrrru]_{i^{F}_{i,m-1}}
}\]
for $i= k$, $k+1$ $\le m-1$,
we have the equation
\[
\{{(e^{F_{m-1}}_{k+1})^{(\ell)}{\times} \{\ast\}}\cup (e^{F_{m-1}}_{k})^{(\ell)}{\times} (e')^{(\ell)}\}{\circ}j_{k+1}{\circ}\hat{\sigma}_{k+1}
  = i^{m,1}_{k+1,m},
\]
where $i^{m,1}_{k+1,m} : F^{m,1}_{k+1}\to F^{m,1}_{m}$ is the inclusion map.
Thus we obtain
\begin{align*}&
 e_{k+1}^{F_m}{\circ}\hat{\lambda}_{k+1}
{\circ}j_{k+1}{\circ}\hat{\sigma}_{k+1}
 = (\mu_{m-1,0}\cup\mu'_{m-1,1}){\circ}i^{m,1}_{k+1,m}
 = i^F_{k+1,m}{\circ}\mu^{m,1}_{k+1}\\
&\qquad = i^F_{k+1,m}{\circ}e_{k+1}^{F_{k+1}}{\circ}\sigma_{k+1}{\circ}\mu^{m,1}_{k+1}
= e_{k+1}^{F_m}{\circ}P^{k+1}\Omega i^F_{k+1,m}{\circ}\sigma_{k+1}{\circ}\mu^{m,1}_{k+1}.
\end{align*}
Hence by the equation \eqref{delta_k+1}, we have
\begin{align*}&
i^F_{k+1,m}{\circ}\mu^{m,1}_{k+1}
  = \nabla_{F_m}{\circ}(e_{k+1}^{F_m}{\circ}\hat{\lambda}_{k+1}{\circ}j_{k+1}{\circ}\hat{\sigma}_{k+1}
\vee e_{k+1}^{F_m}{\circ}\delta_{k+1}){\circ}\nu^{m,1}_{k+1}
\\&\qquad
  = \nabla_{F_m}{\circ}(i^F_{k+1,m}{\circ}\mu^{m,1}_{k+1}
    \vee e_{k+1}^{F_m}\circ\delta_{k+1}){\circ}\nu^{m,1}_{k+1}.
\end{align*}
Using Theorem 2.7 (1) of \cite{Oda92} and 
the multiplication $\mu$ on $G \simeq F_m$,
the map
$e_{k+1}^{F_m}\circ\delta_{k+1} : \Sigma K^{m,1}_{k+1}\to F_m$
is null-homotopic.
Hence by using the exact sequence
\[
[\Sigma K^{m,1}_{k+1}, E^{k+2}\Omega F_m] 
 \xrightarrow{{p_{k+2}^{\Omega F_m}}_{\ast}} [\Sigma K^{m,1}_{k+1}, P^{k+1}\Omega F_m] 
 \xrightarrow{{e^{F_m}_{k+1}}_{\ast}} [\Sigma K^{m,1}_{k+1}, F_m],
\]
we obtain a map $\delta^{\prime}_{k+1} : \Sigma K^{m,1}_{k+1}\to E^{m+1}\Omega F_m$ which satisfies $\delta_{k+1} = p_{k+2}^{\Omega F_m}{\circ}\delta^{\prime}_{k+1}$.
Since
$E^{k+2}\Omega F_m \xrightarrow{p_{k+2}^{\Omega F_m}} P^{k+1}{\Omega F_m}
  \hookxrightarrow{\iota^{\Omega F_m}_{k+1,k+2}} P^{k+2}{\Omega F_m}$
is the cofibre sequence, we have
$\iota^{\Omega F_m}_{k+1,k+2}{\circ}\delta_{k+1}
 = \iota^{\Omega F_m}_{k+1,k+2}{\circ}p_{k+2}^{\Omega F_m}{\circ}\delta^{\prime}_{k+1}
 = 0$
and 
\begin{align*}&
\iota^{\Omega F_m}_{k+1,k+2}{\circ}\nabla_{P^{k+1}\Omega F_m}{\circ}( 
  \hat{\lambda}_{k+1}{\circ}j_{k+1}{\circ}\hat{\sigma}_{k+1} \vee \delta_{k+1}){\circ}\nu^{m,1}_{k+1}
\\&\qquad
 = \nabla_{P^{k+2}\Omega F_m}{\circ}(\iota^{\Omega F_m}_{k+1,k+2} 
  {\circ}\hat{\lambda}_{k+1}{\circ}j_{k+1}{\circ}\hat{\sigma}_{k+1}
     \vee 0){\circ}\nu^{m,1}_{k+1}
\\&\qquad
 = \iota^{\Omega F_m}_{k+1,k+2} 
  {\circ}\hat{\lambda}_{k+1}{\circ}j_{k+1}{\circ}\hat{\sigma}_{k+1}.
\end{align*}
From the equation (\ref{delta_k+1}), we obtain
\[
\iota^{\Omega F_m}_{k+1,k+2}{\circ}P^{k+1}\Omega i^F_{k+1,m}{\circ}\sigma_{k+1}{\circ}\mu^{m,1}_{k+1}
 =\iota^{\Omega F_m}_{k+1,k+2} {\circ}\hat{\lambda}_{k+1}{\circ}j_{k+1}{\circ}\hat{\sigma}_{k+1}.
\]
Therefore we have \eqref{su=us} by induction.
\par
Next, we show the equation
\begin{equation}\label{su=us2}
\iota^{\Omega_{F_m}}_{m,m+1}{\circ}\sigma_m{\circ}\mu^{m,1}_{m}=\iota^{\Omega_{F_m}}_{m,m+1}{\circ}\hat{\lambda}_m{\circ}\hat{\sigma}_m.
\end{equation}
By the condition (1) of Proposition \ref{cd-to-gs} for $\sigma_i$, we obtain 
\[
\sigma_t{\circ}i^F_{t-1,t} = i^{\Omega F_t}_{t-1,t}{\circ}P^{t-1}\Omega i^F_{t-1,t}{\circ}\sigma_{t-1} \qquad \text{for}\quad t=m{-}1, m.
\]
Hence we have
\begin{align*}
\hat{\sigma}_m{\circ}i^{m,1}_{m-1}
& =((\sigma_m{\circ}i^F_{m-1,m}){\times}\{\ast\}\cup (\sigma_{m-1}{\circ}i^F_{m-2,m-1}){\times}\sigma')\\
& = (\iota^{\Omega F_m}_{m-1,m}{\circ}P^{m-1}\Omega i^F_{m-1,m})^{(\ell)}{\circ}\sigma_{m-1}){\times}\{\ast\}\cup \\
& \hspace{30mm} ((\iota^{\Omega F_{m-1}}_{m-2,m-1}{\circ}P^{m-2}\Omega i^F_{m-2,m-1})^{(\ell)}{\circ}\sigma_{m-1}){\times}\sigma'\\
& = \hat {i}_{m-1}{\circ}j_{m-1}{\circ}\hat{\sigma}_{m-1}.
\end{align*}
By the condition (1) of Proposition \ref{cd-to-gs} for $\hat{\lambda}_{i}$, we have $\hat{\lambda}_{m}{\circ}\hat{i}_{m-1}=\iota^{\Omega F_m}_{m-1,m}{\circ}\hat{\lambda}_{m-1}$.
Thus we have the equation
\begin{align*}&
{i^{m,1}_{m-1}}^\ast(\hat{\lambda}_m{\circ}\hat{\sigma}_m)=\hat{\lambda}_m{\circ}\hat{\sigma}_m{\circ}i^{m,1}_{m-1}
= \hat{\lambda}_m{\circ}\hat {i}_{m-1}{\circ}j_{m-1}{\circ}\hat{\sigma}_{m-1}\\
&\qquad= \iota^{\Omega F_m}_{m-1,m}{\circ}\hat{\lambda}_{m-1}{\circ}j_{m-1}{\circ}\hat{\sigma}_{m-1} = \iota^{\Omega F_m}_{m-1,m}{\circ}P^m\Omega i^F_{m-1,m}{\circ}\sigma_{m-1}{\circ}\mu^{m,1}_{m-1}
\\&\qquad = \sigma_m{\circ}i^{F}_{m-1,m}{\circ}\mu^{m,1}_{m} ={i^{m,1}_{m-1}}^\ast(\sigma_m{\circ}\mu^{m,1}_{m})
\end{align*}
by using the equation $k=m{-}1$ of \eqref{su=us}.
Thus by a standard argument of homotopy theory, we obtain the difference map
$\delta_{m} : \Sigma K^{m,1}_{m}\to P^{k+1}\Omega F_m$
such that 
\begin{equation}\label{delta_m}
\sigma_{m}{\circ}\mu^{m,1}_{m} = \nabla_{P^{m}\Omega F_m}
{\circ}(\hat{\lambda}_{m}{\circ}\hat{\sigma}_{m}
\vee \delta_{m}){\circ}\nu^{m,1}_{m}.
\end{equation}
By the condition (2) of Proposition \ref{cd-to-gs} of $\hat{\lambda}_{m}$,
we have the equation
\begin{align*}&
 e_{m}^{F_m}{\circ}\hat{\lambda}_{m}{\circ}\hat{\sigma}_m 
\\&\quad
 =\mu^{m,1}_m{\circ}\{{(e^{F_{m}}_{m})^{(\ell)}{\times} \{\ast\}} \cup (e^{F_{m-1}}_{m-1})^{(\ell)}{\times} (e')^{(\ell)}\}
{\circ}(\sigma_m{\times}\{\ast\}\cup\sigma_{m-1}{\times}\sigma')
 = \mu^{m,1}_m.
\end{align*}
Thus from the equation \eqref{delta_m}, we obtain
\[
\mu^{m,1}_m=\nabla_{F_m}{\circ}( e_{m}^{F_m}{\circ}\hat{\lambda}_{m}{\circ}\hat{\sigma}_m \vee e^{F_m}_m{\circ}\delta_{m}){\circ}\nu^{m,1}_{m}
=\nabla_{F_m}{\circ}( \mu^{m,1}_m \vee e^{F_m}_m{\circ}\delta_{m}){\circ}\nu^{m,1}_{m}.
\]
Thus we obtain $e^{F_m}_m{\circ}\delta_{m}=0$.
By using the fibre sequence $E^{m+1}\Omega F_{m}  \xrightarrow{p_{m+1}^{\Omega F_m}} P^{m}\Omega F_{m} \xrightarrow{e^{F_m}_m} F_{m}$, we obtain a map $\delta'_m:  \Sigma K^{m,1}_{m}\to E^{m+1}\Omega F_{m}$ which satisfies $\delta_m =p^{\Omega F_m}_{m+1}{\circ}\delta_m$.
Since
$E^{m+1}\Omega F_m \xrightarrow{p_{m+1}^{\Omega F_m}} P^{m}{\Omega F_m}
  \hookxrightarrow{\iota^{\Omega F_m}_{m,m+1}} P^{m+1}{\Omega F_m}$
is the cofibre sequence, we have
$\iota^{\Omega F_m}_{m,m+1}{\circ}\delta_{m} = \iota^{\Omega F_m}_{m,m+1}{\circ}p_{m+1}^{\Omega F_m}{\circ}\delta^{\prime}_{m}= 0$.
Thus by the equation (\ref{delta_m}), we obtain
\begin{align*}&
\iota^{\Omega F_m}_{m,m+1}{\circ}\sigma_{m}{\circ}\mu^{m,1}_{m}
= \iota^{\Omega F_m}_{m,m+1}{\circ} \nabla_{P^{m}\Omega F_m}{\circ}(\hat{\lambda}_{m}{\circ}\hat{\sigma}_{m} \vee \delta_{m}){\circ}\nu^{m,1}_{m}
\\&\qquad
=  \nabla_{P^{m+1}\Omega F_m}{\circ}( \iota^{\Omega F_m}_{m,m+1}{\circ} \hat{\lambda}_{m}{\circ}\hat{\sigma}_{m} \vee 0){\circ}\nu^{m,1}_{m}
= \iota^{\Omega F_m}_{m,m+1}{\circ} \hat{\lambda}_{m}{\circ}\hat{\sigma}_{m}.
\end{align*}
This is the proof of \eqref{su=us2}.
\par
Finally,  we construct a map
$\hat{\lambda} : \hat{F}_{m+1} \to P^{m+1}\Omega F_m$.
%
By the condition (1) of Proposition \ref{cd-to-gs} for $\sigma_i$, we obtain \[\sigma_m{\circ}i^F_{m-1,m} = i^{\Omega F_m}_{m-1,m}{\circ}P^{m-1}\Omega i^F_{m-1,m}{\circ}\sigma_{m-1}.\] Hence we have
\begin{align*}
(\sigma_m&{\times}\sigma'){\circ}i^{m,1}_{m}
 = (\sigma_m{\times}\sigma')
  {\circ}(\mathrm{id}_{F_m}{\times}\{\ast\}\cup i^F_{m-1,m}{\times}\mathrm{id}_{F'_{1}})\\
&\hspace{-2mm} = (\mathrm{id}_{P^m\Omega F_m}{\times}\{\ast\}\cup (i^{\Omega F_m}_{m-1,m}{\circ}P^{m-1}\Omega i^F_{m-1,m})^{(\ell)}{\times}\mathrm{id}_{\Sigma\Omega F'_1}){\circ}(\sigma_m{\times}\{\ast\}\cup\sigma_{m-1}{\times}\sigma')\\
&\hspace{-2mm} = \hat{i}_m{\circ}\hat{\sigma}_m.
\end{align*}
Also by condition (1) of Proposition \ref{cd-to-gs} for $\hat{\lambda}_{i}$, we obtain $\hat{\lambda}_{m+1}{\circ}\hat{i}_m = \iota^{\Omega F_m}_{m,m+1}{\circ}\hat{\lambda}_{m}$.
Hence we have 
\[
\hat{\lambda}_{m+1}{\circ}(\sigma_m{\times}\sigma'){\circ}i^{m,1}_{m}
 = \hat{\lambda}_{m+1}{\circ}\hat{i}_m{\circ}\hat{\sigma}_m
 = \iota^{\Omega F_m}_{m,m+1}{\circ}\hat{\lambda}_{m}{\circ}\hat{\sigma}_m.
\]
Moreover, by the equation \eqref{su=us2},
 we have
\begin{align*}
\hat{\lambda}_{m+1}{\circ}(\sigma_m{\times}\sigma'){\circ}i^{m,1}_{m}
 & = \iota^{\Omega F_m}_{m,m+1}{\circ}\sigma_m{\circ}\mu^{m,1}_{m}
 =\iota^{\Omega F_m}_{m,m+1}{\circ}\sigma_m {\circ}\mu'_{m,1}{\circ}i^{m,1}_{m}.
\end{align*}
Hence by using the cofibration $K^{m,1}_{m+1}\xrightarrow{{w}^{m,1}_{m+1}} F^{m,1}_{m} \hookxrightarrow{i^{m,1}_{m}} F^{m,1}_{m+1}$, there is a map 
$\delta_{m+1}: \Sigma K^{m,1}_{m+1}\to P^{m+1}\Omega F_{m}$ such that
\begin{equation}\label{delta_m+1}
\iota^{\Omega F_m}_{m,m+1}{\circ}\sigma_m {\circ}\mu'_{m,1}
= \nabla_{P^{m+1}\Omega F_{m}}{\circ}(\hat{\lambda}_{m+1}{\circ}(\sigma_m{\times}\sigma')
   \vee \delta_{m+1}){\circ}\nu^{m,1}_{m+1}.
\end{equation}
To continue calculating,
we consider the map $\bar{e}: \hat{E}_{m+1} \to \Sigma K^{m+1}_m$
induced from the bottom left square of the following commutative diagram:
\[
\xymatrix{
F^{m,1}_m \ar@{^{(}->}[rr]^{i^{m,1}_m} \ar[d]^{\hat{\sigma}_m}
 &
 & F^{m,1}_{m+1} \ar[rr]^{q} \ar[d]^{\sigma_m {\times} \sigma'}
 &
 & \Sigma K^{m,1}_m \ar[d]^{\Sigma g_m{\ast} g'}
\\
\hat{F}_m \ar@{^{(}->}[rr]^{\hat{i}_m}
    \ar[d]^{\hat{e}_m}
 &
 & \hat{F}_{m+1} \ar[rr]^{\bar{q}} \ar[d]^{(e^{F_m}_m)^{(\ell)} {\times} (e')^{(\ell)}}
 &
 & \hat{E}_{m+1} \ar@{-->}[d]^{\bar{e}}
\\
F^{m,1}_m \ar@{^{(}->}[rr]^{i^{m,1}_m}
 &
 & F^{m,1}_{m+1} \ar[rr]^{q}
 &
 & \Sigma K^{m,1}_m,
}\]
where the map 
$\hat{e}_m : \hat{F}_m \to F^{m,1}_m$ is
$(e^{F_m}_m)^{(\ell)} {\times} \{\ast\}\cup (e^{F_{m-1}}_{m-1})^{(\ell)}{\times} (e')^{(\ell)}$.
Since $\hat{e}_m{\circ}\hat{\sigma}_m$ and 
$((e^{F_m}_m)^{(\ell)} {\times} (e')^{(\ell)}){\circ}(\sigma_m {\times} \sigma')$
are homotopic to the identity maps,
$\bar{e}{\circ}\Sigma g_m{\ast} g_1$ is homotopic to the identity map of $\Sigma K^{m+1}_m$.
Then the equation (\ref{delta_m+1}) is as follows:
\begin{align*}
\iota^{\Omega F_m}_{m,m+1}
 &{\circ}\sigma_m {\circ}\mu'_{m,1}
  = \nabla_{P^{m+1}\Omega F_{m}}{\circ}(\hat{\lambda}_{m+1}{\circ}(\sigma_m{\times}\sigma')
   \vee \delta_{m+1}){\circ}\nu^{m,1}_{m+1}\\
 & = \nabla_{P^{m+1}\Omega F_{m}}{\circ}(\hat{\lambda}_{m+1}{\circ}(\sigma_m{\times}\sigma')
   \vee \delta_{m+1}{\circ}\bar{e}{\circ}\Sigma g_m{\ast} g'){\circ}\nu^{m,1}_{m+1}\\
 & = \nabla_{P^{m+1}\Omega F_{m}}{\circ}(\hat{\lambda}_{m+1}
   \vee \delta_{m+1}{\circ}\bar{e}){\circ}
 ((\sigma_m{\times}\sigma')\vee\Sigma g_m{\ast} g'){\circ}\nu^{m,1}_{m+1}.\\
\intertext{ By Lemma \ref{2cofib}, we proceed further:}
 & = \nabla_{P^{m+1}\Omega F_{m}}{\circ}(\hat{\lambda}_{m+1}
   \vee \delta_{m+1}{\circ}\bar{e}){\circ}\hat{\nu}_{m+1}{\circ}(\sigma_m{\times}\sigma').
\end{align*}
Therefore we adopt
$\nabla_{P^{m+1}\Omega F_{m}}{\circ}(\hat{\lambda}_{m+1}
   \vee \delta_{m+1}{\circ}\bar{e}){\circ}\hat{\nu}_{m+1}$
as
$\hat{\lambda}$.
Then we obtain 
\[
\iota^{\Omega F_m}_{m,m+1}{\circ}\sigma_m {\circ}\mu'_{m,1}=\hat{\lambda}{\circ}(\sigma_m{\times}\sigma').
\]
We show the equation
$e^{F_m}_{m+1}{\circ}\hat{\lambda}
  = \mu'_{m,1}\circ\{(e^{F_m}_m)^{(\ell)}{\times} (e')^{(\ell)}\}$ as follows.
By the condition (2) of \ref{cd-to-gs} of $\hat{\lambda}_{m+1}$, we have
\[
e^{F_m}_{m+1}{\circ}\hat{\lambda}_{m+1}{\circ}(\sigma_m{\times}\sigma')
= \mu'_{m,1}{\circ}\{(e^{F_m}_m)^{(\ell)}{\times} (e')^{(\ell)}\}{\circ}(\sigma_m{\times}\sigma')
= \mu'_{m,1}.
\]
Hence by equations $e^{F_m}_{m+1}{\circ}\iota^{\Omega F_m}_{m,m+1}{\circ}\sigma_m=\mathrm{id}_{F_m}$ and \eqref{delta_m+1}, we have
\begin{align*}
\mu'_{m,1}
&= e^{F_m}_{m+1}{\circ}\nabla_{P^{m+1}\Omega F_{m}}{\circ}(\hat{\lambda}_{m+1}{\circ}(\sigma_m{\times}\sigma')
   \vee \delta_{m+1}){\circ}\nu^{m,1}_{m+1}\\
&= \nabla_{F_m}{\circ}(\mu'_{m,1} \vee e^{F_m}_{m+1} {\circ} \delta_{m+1}){\circ}\nu^{m,1}_{m+1}.
\end{align*}
Thus we obtain the equation $e^{F_m}_{m+1}\circ\delta_{m+1} =0$.
Moreover, we obtain
\begin{align*}
e^{F_m}_{m+1}{\circ}\hat{\lambda}
 & = e^{F_m}_{m+1}{\circ}\nabla_{P^{m+1}\Omega F_{m}}{\circ}(\hat{\lambda}_{m+1}
   \vee \delta_{m+1}{\circ}\bar{e}){\circ}\hat{\nu}_{m+1}\\
 & = \nabla_{F_{m}}{\circ}(e^{F_m}_{m+1}\circ\hat{\lambda}_{m+1}
   \vee 0 ){\circ}\hat{\nu}_{m+1}
 = e^{F_m}_{m+1}\circ\hat{\lambda}_{m+1}
\end{align*}%
and by the condition (2) of Proposition \ref{cd-to-gs} of $\hat{\lambda}_{m+1}$, we obtain
\[
e^{F_m}_{m+1}{\circ}\hat{\lambda}
  = \mu'_{m,1}\circ\{(e^{F_m}_m)^{(\ell)}{\times} (e')^{(\ell)}\}.
\]
This completes the proof.
\end{proof}
Thus we have the following commutative diagram:
\begin{equation}\label{diagPxPmu}
\xymatrix{
F_m \ar[d]^{\iota^{\Omega F_m}_{m,m+1}{\circ}\sigma_m}
 & F_m{\times} A \ar[d]^{\sigma_m{\times} \sigma_A} \ar[l]_{pr_1} \ar[r]^{1{\times} \alpha }
 & F_m{\times} F'_1 \ar[d]^{\sigma_m{\times} \sigma'}  \ar[r]^{\mu'_{m,1}}
 & F_m \ar[d]^{\iota^{\Omega F_m}_{m,m+1}{\circ}\sigma_m}\\
P^{m+1}\Omega F_m \ar[d]^{e^{F_m}_{m+1}}
 & (P^m\Omega F_m)^{(\ell)}{\times} (\Sigma \Omega A)^{(\ell)}
    \ar[l]_{\hspace{-20pt}\phi}
    \ar[d]^{(e^{F_m}_m)^{(\ell)}{\times} (e^A_1)^{(\ell)}}
    \ar[r]^{\hspace{40pt}\chi}
 & \hat{F}_{m+1} \ar[d]^{(e^{F_m}_m)^{(\ell)}{\times} (e')^{(\ell)}} \ar[r]^{\hat{\lambda}}
 & P^{m+1}\Omega F_m \ar[d]^{e^{F_m}_{m+1}}\\
F_m
 & F_m{\times} A \ar[l]_{pr_1} \ar[r]^{1{\times} \alpha }
 & F_m{\times} F'_1 \ar[r]^{\mu'_{m,1}}
 & F_m.
 }\end{equation}
Now we define a cone-decomposition 
$
\{ \hat{E}^{\prime}_{k} \xrightarrow{\hat{w}^{\prime}_{k}}
    \hat{F}^{\prime}_{k-1} \hookxrightarrow{\hat{i}^{\prime}_{k-1}} \hat{F}^{\prime}_{k}
 \,|\, 1 {\le} k {\le} m{+}1 \}
$
of $(P^m\Omega F_m)^{(\ell)}{\times} (\Sigma \Omega A)^{(\ell)}$ of length $m{+}1$ by replacing $F'_1$ with $A$ in the cone-decomposition of
$(P^m\Omega F_m)^{(\ell)}{\times} (\Sigma \Omega F'_1)^{(\ell)}$.

Let us recall the sequence of cofibrations
\begin{equation*}
\{ E^{k}\Omega F_m \xrightarrow{p^{\Omega F_m}_{k}} P^{k-1}\Omega F_m
  \hookxrightarrow{\iota^{\Omega F_m}_{k-1}} P^{k}\Omega F_m
 \,|\, 1\le k \le m{+}1 \}
\end{equation*}
which gives a cone-decomposition of $P^{m+1}\Omega F_m$ of length $m+1$.
Let $D$ be the homotopy pushout of 
$(\iota^{\Omega F_m}_{m,m+1})^{(\ell)}{\circ}pr_1$ and 
$\hat{\lambda}{\circ}(\mathrm{id}_{(P^m\Omega F_m)^{(\ell)}}
 {\times} (\Sigma \Omega \alpha)^{(\ell)})$:
\[
\xymatrix{
(P^m\Omega F_m)^{(\ell)}{\times} (\Sigma \Omega A)^{(\ell)}
 \ar[d]^{f^{\gets}}
 \ar[r]^{\hspace{30pt}f^{\to}}
 & P^{m+1}\Omega F_m \ar[d]\\
P^{m+1}\Omega F_m \ar[r]
 & D.
}\]
Here $f^{\to}
 = \hat{\lambda}{\circ}(\mathrm{id}_{(P^m\Omega F_m)^{(\ell)}} {\times} (\Sigma \Omega \alpha)^{(\ell)})$
and $f^{\gets} = (\iota^{\Omega F_m}_{m,m+1})^{(\ell)}{\circ}pr_1$.
We define a cone-decomposition of $D$ as follows.
By the equation
$
\hat{\lambda}{\circ}\hat{i}_{m}
  = \nabla_{P^{m+1}\Omega F_{m}}{\circ}(\hat{\lambda}_{m+1}
   \vee \delta_{m+1}{\circ}\bar{e}){\circ}\hat{\nu}_{m+1}{\circ}\hat{i}_{m}
  = \hat{\lambda}_{m+1}{\circ}\hat{i}_{m},
$
we do not distinguish the restriction of $\hat{\lambda}$ on $\hat{F}_{k}$ and $\hat{\lambda}_{k}$ and hence $f^{\to}$ is a filtered map up to homotopy.
Since 
$\hat{E}^{\prime}_{k} \xrightarrow{\hat{w}^{\prime}_{k}}
    \hat{F}^{\prime}_{k-1} \hookxrightarrow{\hat{i}^{\prime}_{k-1}} \hat{F}^{\prime}_{k}$ is a cofibre sequence, we have
\begin{align*}
e^{F_m}_{k-1}{\circ}(f^{\to}|_{\hat{F}^{\prime}_{k-1}}{\circ}\hat{w}^{\prime}_{k})
 & = e^{F_m}_{k}{\circ}\iota^{\Omega F_m}_{k-1,k}
    {\circ}f^{\to}|_{\hat{F}^{\prime}_{k-1}}{\circ}\hat{w}^{\prime}_{k}\\
 & =  e^{F_m}_{k}{\circ}f^{\to}|_{\hat{F}^{\prime}_{k}}
   {\circ}\hat{i}^{\prime}_{k-1}{\circ}\hat{w}^{\prime}_{k}
  =  e^{F_m}_{k}{\circ}f^{\to}|_{\hat{F}^{\prime}_{k}}{\circ}0  = 0.
\end{align*}
Using the fibre sequence 
$E^{k}\Omega F_m \xrightarrow{p^{\Omega F_m}_{k}} P^{k-1}\Omega F_m
 \xrightarrow{e^{F_m}_{k-1}} F_m$,
there exists a map $g^{\to}_k : \hat{E}^{\prime}_{k} \to E^{k}\Omega F_m$
such that the commutativity of the following diagram:
\begin{equation}\label{g^right}
\xymatrix{
\hat{E}^{\prime}_{k} \ar[rr]^{\hat{w}^{\prime}_{k}} \ar[d]^{g^{\to}_k}
 && \hat{F}^{\prime}_{k-1} \ar@{^{(}->}[rr]^{\hat{i}^{\prime}_{k-1}}
   \ar[d]^{f^{\to}|_{\hat{F}^{\prime}_{k-1}}}
 && \hat{F}^{\prime}_{k}\ar[d]^{f^{\to}|_{\hat{F}^{\prime}_{k}}}\\
E^{k}\Omega F_m \ar[rr]^{p^{\Omega F_m}_{k}}
 && P^{k-1}\Omega F_m \ar@{^{(}->}[rr]^{\iota^{\Omega F_m}_{k-1,k}}
 && P^{k}\Omega F_m.
}
\end{equation}
Since $f^{\gets}$ is composition of the projection and the inclusion,
it is clear that there exists a map
$g^{\gets}_k : \hat{E}^{\prime}_{k} \to E^{k}\Omega F_m$
satisfy that the following diagram is commutative:
\begin{equation}\label{g^left}
\xymatrix{
\hat{E}^{\prime}_{k} \ar[rr]^{\hat{w}^{\prime}_{k}} \ar[d]^{g^{\gets}_k}
 && \hat{F}^{\prime}_{k-1} \ar@{^{(}->}[rr]^{\hat{i}^{\prime}_{k-1}}
   \ar[d]^{f^{\gets}|_{\hat{F}^{\prime}_{k-1}}}
 && \hat{F}^{\prime}_{k}\ar[d]^{f^{\gets}|_{\hat{F}^{\prime}_{k}}}\\
E^{k}\Omega F_m \ar[rr]^{p^{\Omega F_m}_{k}}
 && P^{k-1}\Omega F_m \ar@{^{(}->}[rr]^{\iota^{\Omega F_m}_{k-1,k}}
 && P^{k}\Omega F_m.
}
\end{equation}
Let $E^{P}_{k}$ be a homotopy pushout of $g^{\to}_k$ and $g^{\gets}_k$,
 and $F^{P}_{k}$ be a homotopy pushout of $f^{\to}|_{\hat{F}^{\prime}_{k}}$
and $f^{\gets}|_{\hat{F}^{\prime}_{k}}$,
then using diagrams (\ref{g^right}) and (\ref{g^left}) and using the universal property of the homotopy pushout, we have the following diagram such that the front column
$ E^{P}_{k} \to F^{P}_{k-1} \to F^{P}_{k}$ is a cofibration:
\[
\xymatrix{
& & \hat{E}^{\prime}_k
  \ar[lld]_(.2){g^{\gets}_k}
  \ar[d]^{\hat{w}^{\prime}_{k}}
  \ar[rrrd]^{g^{\to}_k}
& & & 
\\
E^{k}\Omega F_m
  \ar[d]^{p^{\Omega F_m}_{k}}
  \ar[rrrd]
& &  \hat{F}^{\prime}_{k-1}
  \ar[lld]|(.4)\hole _(.2){f^{\gets}|_{\hat{F}^{\prime}_{k-1}}}
  \ar@{^{(}->}[d]|(.63)\hole ^(.4){\hat{i}^{\prime}_{k-1}}
  \ar[rrrd]|(.61)\hole ^{f^{\to}|_{\hat{F}^{\prime}_{k-1}}}
& & & E^{k}\Omega F_m
  \ar[lld]
  \ar[d]^{p^{\Omega F_m}_{k}}
\\
P^{k-1}\Omega F_m
  \ar@{^{(}->}[d]^{\iota^{\Omega F_m}_{k-1,k}}
  \ar[rrrd]
& & \hat{F}^{\prime}_k
  \ar[lld]|(.4)\hole _(.2){f^{\gets}|_{\hat{F}^{\prime}_{k}}}
  \ar[rrrd]|(.36)\hole |(.61)\hole ^{f^{\to}|_{\hat{F}^{\prime}_{k}}}
 & E^{P}_{k}
  \ar[d]
& & P^{k-1}\Omega F_m
  \ar[lld]
  \ar@{^{(}->}[d]^{\iota^{\Omega F_m}_{k-1,k}}
\\
P^{k}\Omega F_m
  \ar[rrrd]
& & &  F^{P}_{k-1}
  \ar@{^{(}->}[d]
& & P^{k-1}\Omega F_m
  \ar[lld]
\\
& &
&  F^{P}_{k}.
& &
\\
}
\]
Thus we obtain
a cone-decomposition of $D$ of length $m+1$: 
\begin{equation*}
\{ E^{P}_{k} \to F^{P}_{k-1} \hookrightarrow F^{P}_{k}
\,|\, 1\le k \le m+1\}.
\end{equation*}
Therefore we have the inequalities
\begin{equation*}\label{cat-cl-D}
\cat{D}\le\cl{D}\le m+1.
\end{equation*}
Recall the horizontal top and bottom lines of the diagram (\ref{diagPxPmu}).
The homotopy pushout of these lines are $G\cup_{\psi}G{\times} CA$.
Since dimensions of $F_m$, $F_1$ and $A$ are less than or equal to $\ell$,
all composition of columns in the diagram (\ref{diagPxPmu}) are homotopic to identity maps.
By the universal property of the homotopy pushout,
we obtain a composite map
$D\to G\cup_{\psi}G{\times} CA \simeq E \to D$
which is homotopic to the identity map.
Thus $D$ dominates $E$ and we have
\[
\cat{E}\le\cat{D}\le\cl{D}\le m+1.
\]
This completes the proof of Theorem \ref{main2thm}.
\section{Application of Theorem \ref{main2thm}}%
\label{appli-mthm}
We want to determine the L-S category of $\mathrm{SO}(10)$
by applying the principal bundle $p : \mathrm{SO}(10)\to {S}^9$ to
Theorem \ref{main2thm}.
First, we estimate the lower bound of cat(SO(10)). 
For $k=\field_{2}$ the prime field of characteristic 2, the cohomology ring of $\mathrm{SO}(10)$ with coefficients in $k$ is well-known and described as follows:
\[
H^*(\mathrm{SO}(10);k)\cong P_k[x_1, x_3, x_5, x_7, x_9 ]/(x_1^{16}, x_3^4, x_5^2, x_7^2, x_9^2),
\]
where the suffix of the generator indicates its dimension. Hence, we have
\begin{equation}\label{eq:lower-bound-so10}
21 = \cupln{\mathrm{SO}(10)}{k} \leq \cat{\mathrm{SO}(10)}. 
\end{equation}
\par
To estimate the upper bound by using Theorem \ref{main2thm}, let us recall the cone-decomposition of $\mathrm{SO}(9)$: the cone-decomposition of $\mathrm{Spin}(7)$ is given by Iwase, Mimura and Nishimoto\cite{IMN03} as follows:
\[
\ast \subset F^{\prime}_1=\Sigma \mathbb{C}\mathrm{P}^3
     \subset F^{\prime}_2\subset F^{\prime}_3
     \subset F^{\prime}_4\subset F^{\prime}_5\simeq \mathrm{Spin}(7) .
\]
In Iwase, Mimura and Nishimoto \cite{IMN05}, the cone-decomposition of $\mathrm{SO}(9)$ 
\[
\{K_{i}\to F_{i-1}\to F_{i} \;\;|\;\; 1\le i \le 20,\; F_{0}= \{\ast\}
   \text{\; and\; }F_{20}= \mathrm{SO}(9) \}
\]
is given by using the filtration ${F'_i}$ and  principal bundle 
$\mathrm{Spin}(7)\hookrightarrow \mathrm{SO}(9)\to \mathbb{R}\mathrm{P}^{15}$.
Let $e^k$ be the $k$-dimension cell of $\mathbb{R}\mathrm{P}^{15}$, then $F_i$ can be described as
\begin{equation*}
\begin{aligned}
F_1 & = F'_1 \cup e^1 =\Sigma\mathbb{C}\mathrm{P}^3\vee \mathrm{S}^1  \\
F_2 & = F'_2 \cup (e^1 {\times} F'_1) \cup e^2 \\
F_3 & = F'_3 \cup (e^1 {\times} F'_2) \cup (e^2 {\times} F'_1) \cup e^3 \\
F_4 & = F'_4 \cup (e^1 {\times} F'_3) \cup (e^2 {\times} F'_2) \cup (e^3 {\times} F'_1) \cup e^4 \\
F_5 & = F'_5 \cup (e^1 {\times} F'_4) \cup (e^2 {\times} F'_3) \cup (e^3 {\times} F'_2) \cup (e^4 {\times} F'_1) \cup e^5 \\
F_6 & = F'_5 \cup (e^1 {\times} F'_5) \cup (e^2 {\times} F'_4) \cup (e^3{\times} F'_3) \cup (e^4{\times} F'_2) \cup (e^5{\times} F'_1) \cup e^6 \\
F_7 & = F'_5 \cup (e^1 {\times} F'_5) \cup (e^2 {\times} F'_5) \cup (e^3 {\times} F'_4) 
\cup (e^4{\times} F'_3) \cup (e^5{\times} F'_2) \cup (e^6{\times} F'_1) \cup e^7 \\
 \vdots & \hspace*{30mm}\vdots  \\
F_{15} & = F'_5 \cup (e^1 {\times} F'_5) \cup\dots\cup (e^{10} {\times} F'_5) \cup (e^{11} {\times} F'_4) \cup\dots\cup (e^{14}{\times} F'_1) \cup e^{15} \\
F_{16} & = F'_5 \cup (e^1 {\times} F'_5) \cup\dots\cup (e^{11} {\times} F'_5) \cup (e^{12} {\times} F'_4) \cup\dots\cup (e^{15}{\times} F'_1) \\
F_{17} & = F'_5 \cup (e^1 {\times} F'_5) \cup\dots\cup (e^{12} {\times} F'_5) \cup (e^{13} {\times} F'_4) \cup (e^{14} {\times} F'_3) \cup (e^{15} {\times} F'_2) \\
F_{18} & = F'_5 \cup (e^1 {\times} F'_5) \cup\dots\cup (e^{13} {\times} F'_5) \cup (e^{14} {\times} F'_4) \cup (e^{15} {\times} F'_3) \\
F_{19} & = F'_5 \cup (e^1 {\times} F'_5) \cup\dots\cup (e^{14} {\times} F'_5) \cup (e^{15} {\times} F'_4) \\
F_{20} & = F'_5 \cup (e^1 {\times} F'_5) \cup\dots\cup (e^{15} {\times} F'_5).
\end{aligned}
\end{equation*}
By the fact $F'_5=F'_4\cup e^{21}$, we have $K_{20}=S^{35}$, where $e^{21}$ is the top cell of $\mathrm{Spin}(7)$.
Since $\mu\vert_{F'_{i}{\times}F'_{1}}$ is compressible into $F'_{i+1}$ for $1 \leq i < 5$ from the proof of Theorem 2.9 of  \cite{IMN05}, $\mu'\vert_{F_{i}{\times}F'_{1}}$ is compressible into $F_{i+1}$ for  $1 \leq i < 20$, 
where $\mu$ and $\mu'$ are multiplications of $\mathrm{Spin}(7)$ and $\mathrm{SO}(9)$, respectively.
Let us consider principal bundles $p : \mathrm{SO}(10)\rightarrow S^9$ and $p^{\prime} : \mathrm{SU}(5)\rightarrow S^9$, together with the following commutative diagram:
\[\xymatrix{
\Sigma \mathbb{C}\mathrm{P}^3 \ar@{^{(}->}[r]
 & \mathrm{SU}(4) \ar@{^{(}->}[r] \ar@{^{(}->}[d]
 & \mathrm{SO}(9) \ar@{^{(}->}[d]
\\
 & \mathrm{SU}(5) \ar@{^{(}->}[r] \ar[rd]^{p^{\prime}}
 & \mathrm{SO}(10) \ar[d]^{p}
\\
 & S^8 \ar@{^{(}->}[r] \ar[uur]|(.31)\hole|\hole ^(.7){\alpha}
       \ar@/^2pc/@{-->}[uu]^{\alpha'} \ar@/^1pc/@{-->}[luu]^{\Sigma{\gamma_{3}}}
 & S^9.
}\]
Here $\alpha: S^8 \to \mathrm{SO}(9)$ is a characteristic map of the principal bundle $p : \mathrm{SO}(10)\rightarrow S^9$.
By Steenrod \cite{Steenrod51}, $\alpha$ is homotopic in $\mathrm{SO}(9)$ to $\alpha^{\prime} : S^8 \rightarrow \mathrm{SU}(4)$ the characteristic map of the principal bundle $p' : \mathrm{SU}(5) \to S^{9}$.
Further by Yokota \cite{Yokota71}, the suspension $\Sigma \gamma_3 : S^8 \to \Sigma \mathbb{C}\mathrm{P}^3$ of the canonical $S^{1}$ bundle $\gamma_{3} : S^{7} \to \mathbb{C}\mathrm{P}^3$
gives a representative of the homotopy class $\alpha^{\prime}$.
Therefore the characteristic map $\alpha$ is compressible 
into $\Sigma \mathbb{C}\mathrm{P}^3\subset F_1$.
Since $\alpha$ is given by a suspension map, we obtain $H_1(\alpha )=0
\in \pi_{8}(\Omega\Sigma\mathbb{C}\mathrm{P}^3{\ast}\Omega\Sigma\mathbb{C}\mathrm{P}^3)$.
Hence by Theorem \ref{main2thm} with $F'_{1}=\Sigma\mathbb{C}\mathrm{P}^3$, we obtain 
\begin{equation}\label{eq:upper-bound-so10}
\cat{\mathrm{SO}(10)} \leq 20 + 1=21. 
\end{equation}
Combining (\ref{eq:upper-bound-so10}) with (\ref{eq:lower-bound-so10}), we obtain the following result.
\begin{thm}\label{catso10=21}
$\cat{\mathrm{SO}(10)}=21$.
\end{thm}

%

\end{document}